\documentclass[leqno,12pt]{article} 

\usepackage{amsmath}
\usepackage{amssymb}
\usepackage{pst-all}

\newtheorem{theorem}{Theorem}[section]

\newtheorem{proposition}[theorem]{Proposition}
\newtheorem{corollary}[theorem]{Corollary}

\newtheorem{lemma}[theorem]{Lemma}

\newcommand{\C}{\mathbb{C}}
\newcommand{\R}{\mathbb{R}}
\newcommand{\N}{\mathbb{N}}

\newcommand{\qed}{\ensuremath{\quad \blacksquare} \vspace{\baselineskip}}

\newcommand{\Diag}{\mathrm{Diag\,}}
\newcommand{\diag}{\mathrm{diag\,}}

\begin{document}
\allowdisplaybreaks
\title{\textbf{On the Higher-Order Derivatives of Spectral Functions: Two Special Cases}}
\author{\textbf{
Hristo S.\ Sendov\thanks{Department of Mathematics and Statistics,
University of Guelph, Guelph, Ontario, Canada N1G 2W1. Email:
\texttt{hssendov\char64 uoguelph.ca}. Research supported by
NSERC.}
 }}

\maketitle

\begin{abstract}
In this work we use the tensorial language developed in
\cite{Sendov2003a} and \cite{Sendov2003b} to differentiate
functions of eigenvalues of symmetric matrices. We describe the
formulae for the $k$-th derivative of such functions in two cases.
The first case concerns the derivatives of the composition of an
arbitrary differentiable function with the eigenvalues at a matrix
with distinct eigenvalues. The second development describes the
derivatives of the composition of a separable symmetric function
with the eigenvalues at an arbitrary symmetric matrix. In the
concluding section we re-derive the formula for the Hessian of a
general spectral function at an arbitrary point. Our approach
leads to a shorter, streamlined derivation than the original in
\cite{LewisSendov:2000a}. The language we use, based on the
generalized Hadamard product, allows us to view the
differentiation of spectral functions as a routine calculus-type
procedure.
\end{abstract}



\section{Introduction}

We say that a real-valued function $F$, on a symmetric matrix
argument, is spectral if it has the following invariance property:
$$
F(UXU^T)=F(X),
$$
for every symmetric matrix $X$ in its domain and every orthogonal
matrix $U$.  The restriction of $F$ to the subspace of diagonal
matrices defines (almost) a function $f(x):=F(\Diag x)$ on a
vector argument $x \in \R^n$. It is easy to see that $f: \R^n
\rightarrow \R$ has the property
$$
f(x)=f(Px) \mbox{  for any permutation matrix  }  P \mbox{ and any
} x \in \mbox{domain} \, f.
$$
We call such functions {\it symmetric}. It is not difficult to see
that $F(X) = f(\lambda(X))$, where $\lambda(X)$ is the vector of
eigenvalues of $X$. An important subclass of spectral functions is
obtained when $f(x) = g(x_1) + \cdots + g(x_n)$ for some function
$g$ on one real variable. We call such symmetric functions {\it
separable} and their corresponding spectral functions will be
called {\it separable spectral} functions.

In \cite{Lewis:1994b} an explicit formulae for the gradient the
spectral function $F$ in terms of the derivatives of the symmetric
function $f$ was given:
\begin{align} \label{deriv:1} \nabla (f \circ \lambda)(X)
&= V \big(\Diag \nabla f(\lambda(X)) \big) V^T,
\end{align}
where $V$ is any orthogonal matrix such that $X=V \big(\mbox{\rm
Diag}\, \lambda(X) \big) V^T$ is the {\it ordered} spectral
decomposition of $X$. In \cite{LewisSendov:2000a} a formula for
the Hessian of $F$ was given, whose structure appeared quite
different than the one for the gradient.  In this work we
generalize the work in \cite{Lewis:1994b} and
\cite{LewisSendov:2000a} by proving the following formula for the
$k$-th derivative of a spectral function
\begin{equation}
\label{mnfrmla} \nabla^k (f \circ \lambda)(X) = V\Big(
\sum_{\sigma \in P^{k}} \Diag^{\sigma}
\mathcal{A}_{\sigma}(\lambda(X)) \Big)V^T,
\end{equation}
 where
again $X=V(\Diag \lambda(X))V^T$. The sum is taken over all
permutations on $k$ elements. (The role of the permutations is
just as a convenient tool for enumerating the maps
$\mathcal{A}_{\sigma}(x)$.) The precise meaning of the operators
$\Diag^{\sigma}$, generalizing the $\Diag$ operator, is explained
in the next section, see Formula~(\ref{diagonal}). The main thing
to keep in mind about the formula is that the maps
$\mathcal{A}_{\sigma}(x)$ depend only on the partial derivatives
of $f(x)$, up to order $k$, and do not depend on the eigenvalues.
In this sense the process of differentiating $f \circ \lambda$
leaves the eigenvalues unscathed, since the only way in which they
participate in the formula above is through the compositions
$\mathcal{A}_{\sigma}(\lambda(X))$ and the conjugation by the
orthogonal matrix $V$.

We show that Formula~(\ref{mnfrmla}) holds in two general
subcases. It holds when $f$ is a $k$-times differentiable
function, not necessarily symmetric, and $X$ is a matrix with
distinct eigenvalues. It also holds when $f$ is a separable
symmetric function and $X$ is an arbitrary symmetric matrix. We
give an easy recipe for computing the maps
$\mathcal{A}_{\sigma}(x)$ in the above two cases.

In addition, we show that in the case when $f$ is a $k$-times {\it
continuously} differentiable, separable, symmetric function,
Formula~(\ref{mnfrmla}) can be significantly simplified. In that
case, all the maps $\mathcal{A}_{\sigma}(x)$ coincide, that is
$\mathcal{A}_{\sigma_1}(x) = \mathcal{A}_{\sigma_2}(x)$ for any
two permutations $\sigma_1, \sigma_2$ on $k$ elements.

Finally, in the last section, we re-derive the formula for the
Hessian of a general spectral function at an arbitrary symmetric
matrix. Our approach leads to a shorter, more streamlined
derivation than the original derivation in
\cite{LewisSendov:2000a}.

The language that we use, based on the generalized Hadamard
product, allows us to differentiate Formula~(\ref{mnfrmla}) just
like one would expect: writing the differential quotient and
taking the limit as the perturbation goes to zero. This gives a
clear view of where the different pieces in the differential come
from and give the process a routine Calculus-like flavour.

In the next section, we give all the necessary notation,
definitions, and background results to make the reading of this
work self-contained.

\section{Notation and background results}
\label{NABR}

We use pretty much the same notation as in the preceding two
papers \cite{Sendov2003a} and \cite{Sendov2003b}. We will briefly
summarize it here for completeness and will try to make the
reading of this part independent.

By $S^n$, $O^n$, and $P^n$ we will denote the set of all $n \times
n$ real symmetric, orthogonal, permutation matrices respectively.
By $M^n$ will be denoted the real Euclidean space of all $n \times
n$ matrices with inner product $\langle X,Y \rangle =
\mbox{tr}\,(XY^T)$. For $A \in S^n$, $\lambda(A)=(\lambda_1(A),
...,\lambda_n(A))$ will be the vector of its eigenvalues ordered
in nonincreasing order. By $\N_k$ we will denote the set
$\{1,2,...,k\}$. For any vector $x$ in $\R^n$, Diag$\,x$ will
denote the diagonal matrix with the vector $x$ on the main
diagonal, and diag$\,: M^n \rightarrow \R^n$ will denote its
conjugate operator, defined by diag$\,(X)=(x_{11},...,x_{nn})$. By
$\R_{\downarrow}^n$ we denote the cone of all vectors $x$ in
$\R^n$ such that $x_1 \ge x_2 \ge \cdots \ge x_n$.
Denote the standard basis in $\R^n$ by $e^1,e^2,...,e^n$. For a
permutation matrix $P \in P^n$ we say that $\sigma : \N_n
\rightarrow \N_n$ is its corresponding permutation map and write
$P \leftrightarrow \sigma$ if for any $h \in \R^n$ we have
$Ph=(h_{\sigma(1)},...,h_{\sigma(n)})^T$ or, in other words,
$P^Te^i = e^{\sigma(i)}$ for all $i=1,...,n$. The symbol
$\delta_{ij}$ will denote the Kroneker delta. It is equal to one
if $i=j$ and zero otherwise.

Any vector $\mu \in \R^n$ defines a partition of $\N_n$ into
disjoint {\it blocks}, where integers $i$ and $j$ are in the same
block if, and only if, $\mu_i=\mu_j$. In general, the blocks that
$\mu$ determines need not contain consecutive integers. We agree
that the block containing the integer $1$ will be the first block,
$I_1$, the block containing the smallest integer that is not in
$I_1$ will be the second block, $I_2$, and so on. The number $r$
will denote the number of blocks in the partition. Let $\iota_l$
denote the largest
 integer in $I_l$ for all $l = 1,...,r$. For any two integers, $i,j \in \N_n$
we will say that they are {\it equivalent (with respect to $\mu$)}
and write $i \sim j$ (or $i \sim_{\mu} j$) if $\mu_i = \mu_j$,
that is, if they are in the same block. Two $k$-indexes
$(i_1,...,i_k)$ and $(j_1,...,j_k)$ are called {\it equivalent} if
$i_l \sim j_l$ for all $l=1,2,...,n$, and we will write
$$
(i_1,...,i_k) \sim (j_1,...,j_k).
$$
A {\it $k$-tensor, $T$, on $\R^n$} is a map from $\R^n \times
\cdots \times \R^n$ ($k$-times) to $\R$ that is linear in each
argument separately. Denote the set of all $k$-tensors on $\R^n$
by $T^{k,n}$. The value of the $k$-tensor at $(h_1,...,h_k)$ will
be denoted by $T[h_1,...,h_k]$. The tensor is called {\it
symmetric} if for any permutation, $\sigma$, on $\N_k$ it
satisfies $T[h_{\sigma(1)},...,h_{\sigma(k)}] = T[h_1,...,h_k]$,
for any $h_1,...,h_k \in \R^n$. Given a vector $\mu \in \R^n$, a
tensor $T \in T^{k,n}$ is {\it $\mu$-symmetric} if for any
permutation $P \in P^n$, such that $P\mu = \mu$, we have
$T[Ph_{1},...,Ph_{k}] = T[h_1,...,h_k], \mbox{ for any }
h_1,...,h_k \in \R^n$. A $k$-tensor valued map, $\mu \in \R^n
\rightarrow \mathcal{F}(\mu) \in T^{k,n}$, is {\it
$\mu$-symmetric} if for every $\mu \in \R^n$ and permutation
matrix $P$ we have $\mathcal{F}(P\mu)[Ph_1,...,Ph_k] =
\mathcal{F}(\mu)[h_{1},...,h_{k}], \mbox{ for any } h_1,...,h_k
\in \R^n$. The tensor is called is called {\it
$\mu$-block-constant} if $T^{i_1...i_k}=T^{j_1...j_k}$ whenever
$(i_1,...,i_k) \sim (j_1,...,j_k)$. A $k$-tensor valued map, $\mu
\in \R^n \rightarrow \mathcal{F}(\mu) \in T^{k,n}$, is {\it
block-constant} if $\mathcal{F}(\mu)$ is $\mu$-block-constant for
every $\mu$. Clearly, every $\mu$-block-constant tensor is
$\mu$-symmetric.  By $T[h]$ we denote the $(k-1)$-tensor on $\R^n$
given by $T[\cdot,...,\cdot,h]$. Similarly for $T[M]$, if $T$ is a
$k$-tensor on $M^n$ and $M \in M^n$.  The following easy lemma was
proved in \cite{Sendov2003b}.
\begin{lemma}
\label{sec1-lastlem} If a $k$-tensor valued map, $\mu \in \R^n
\rightarrow T(\mu) \in T^{k,n}$, is $\mu$-symmetric and
differentiable, then its differential, $\nabla T(\mu)$, is also
$\mu$-symmetric.
\end{lemma}

For each permutation $\sigma$ on $\N_k$, we define {\it
$\sigma$-Hadamard} product between $k$ matrices to be a $k$-tensor
on $\R^n$ as follows. Given any $k$ basic matrices $H_{p_1q_1}$,
$H_{p_2q_2}$,...,$H_{p_kq_k}$
\begin{align*}
 (H_{p_1q_1} \circ_{\sigma} H_{p_2q_2}
\circ_{\sigma} \cdots \circ_{\sigma} H_{p_kq_k})^{i_1i_2...i_k} &=
\left\{
\begin{array}{ll}
1, \hspace{-0.2cm} & \mbox{ if } i_s=p_s=q_{\sigma(s)}, \forall s=1,...,k, \\
0, \hspace{-0.2cm} & \mbox{ otherwise. }
\end{array}
\right.
\end{align*}

Extend this product to a multi-linear map on $k$ matrix arguments:
\begin{align}
\label{ext-by-lin} (H_{1}\circ_{\sigma} H_{2} \circ_{\sigma}
\cdots \circ_{\sigma} H_k)^{i_1i_2...i_k} &=
H_1^{i_1i_{\sigma^{-1}(1)}}\cdots H_k^{i_ki_{\sigma^{-1}(k)}}.
\end{align}
Notice that when $k=1$ we have
$\circ_{\hspace{-0.05cm}_{\mbox{\tiny $(1)$}}} H = \diag H$. Let
$T$ be an arbitrary $k$-tensor on $\R^n$ and let $\sigma$ be a
permutation on $\N_k$. We define $\Diag^{\sigma} T$ to be a
$2k$-tensor on $\R^n$ in the following way
\begin{align}
\label{diagonal} (\Diag^{\sigma} T)^{\substack{i_1...i_k \\
j_1...j_k}} = \left\{
\begin{array}{ll}
T^{i_1...i_k}, & \mbox{ if } i_s = j_{\sigma(s)}, \forall s = 1,...,k, \\
0, & \mbox{ otherwise. }
\end{array}
\right.
\end{align}
When $k=1$ we have $\Diag^{\mbox{\tiny $(1)$}} x = \Diag x$ for
any $x \in \R^n$. Any $2k$-tensor, $T$, on $\R^n$ can naturally be
viewed as a $k$-tensor on $M^n$ in the following way
$$
T[H_1,...,H_k]= \sum_{p_1,q_1=1}^n \cdots \sum_{p_k,q_k=1}^n
T^{\substack{p_1...p_k \\ q_1...q_k}} H_1^{p_1q_1} \cdots
H_k^{p_kq_k}.
$$
Define dot product between two tensors in $T^{k,n}$ in the usual
way:
$$
\langle T_1, T_2\rangle = \sum_{p_1,...,p_k = 1}^n
T_1^{p_1...p_k}T_2^{p_1...p_k}.
$$
 We define an action (called {\it conjugation}) of the orthogonal group $O^n$ on the
space of all $k$-tensors on $\R^n$. For any $k$-tensor, $T$, and
$U \in O^n$ this action will be denoted by $UTU^T \in T^{k,n}$:
\begin{equation}
\label{kur0} (UTU^T)^{i_1...i_k} = \sum_{p_1 = 1 }^{n} \cdots
\sum_{p_k =1}^{n} \Big( T^{p_1...p_k} U^{i_{1}p_{1}} \cdots
U^{i_{k}p_{k}} \Big).
\end{equation}
We showed in \cite{Sendov2003a} that this action is norm
preserving and associative: $V(UTU^T)V^T = (VU)T(VU)^T$ for all
$U,V \in O^n$.

The $\Diag^{\sigma}$ operator, the $\sigma$-Hadamard product, and
conjugation by an orthogonal matrix are connected by the following
formula, see \cite{Sendov2003a}.
\begin{theorem}
\label{jen3a} For any $k$-tensor $T$, any matrices
$H_1$,...,$H_k$, any orthogonal matrix $V$, and any permutation
$\sigma$ in $P^k$ we have the identity
\begin{equation}
\label{gen:1} \langle T, \tilde{H}_1 \circ_{\sigma} \cdots
\circ_{\sigma} \tilde{H}_k \rangle = \big(V(\Diag^{\sigma}
T)V^T\big)[H_1,...,H_k],
\end{equation}
where $\tilde{H}_i = V^TH_iV$, $i=1,...,k$.
\end{theorem}
We will also need the following lemma from \cite{Sendov2003a}.
\begin{lemma}
\label{simplelem} Let $T$ be any $2k$-tensor on $R^n$, $U \in
O^n$, and let $H$ be any matrix. Then, the following identity
holds.
$$
U(T[U^THU])U^T = (UTU^T)[H].
$$
\end{lemma}

Given a permutation $\sigma$ on $\N_k$ we can naturally view it as
a permutation on $\N_{k+1}$ fixing the last element. Let $\tau_l$
 be the transposition $(l,k+1)$, for all $l=1,...,k,k+1$.
 Define $k+1$ permutations, $\sigma_{(l)}$, on $\N_{k+1}$, as follows:
\begin{equation}
\label{permdefn} \sigma_{(l)} = \sigma \tau_l, \mbox{ for }
l=1,...,k,k+1.
\end{equation}
Informally speaking, given the cycle decomposition of $\sigma$, we
obtain $\sigma_{(l)}$, for each $l=1,...,k$, by inserting the
element $k+1$ immediately after the element $l$, and when $l=k+1$,
the permutation $\sigma_{(k+1)}$ fixes the element $k+1$. Clearly
$\sigma_{(l)}^{-1}(k+1)=l$ for all $l$, and
$$
\{ \mbox{All permutations on } \N_{k+1} \} = \{\sigma \tau_l \, |
\, \sigma \mbox{ is a permutation on } \N_k, \,\,\,
l=1,...,k,k+1\}.
$$
For a fixed vector $\mu \in \R^n$ we define $k$ linear maps
$$
T \in T^{k,n} \rightarrow T^{(l)}_{\mbox{\rm \scriptsize out}} \in
T^{k+1,n}, \mbox{ for } l=1,2,...,k,
$$
as follows:
\begin{equation}
\label{defn-out} \big(T^{(l)}_{\mbox{\rm \scriptsize
out}}\big)^{i_1...i_k i_{k+1}} = \left\{
\begin{array}{ll}
0, & \mbox{ if } i_l \sim i_{k+1}, \\
\displaystyle\frac{T^{i_1...i_{l-1}i_{k+1}i_{l+1}...i_k} -
T^{i_1...i_{l-1}i_li_{l+1}...i_k}}{\mu_{i_{k+1}} - \mu_{i_l}}, &
\mbox{ if } i_l \not\sim i_{k+1}.
\end{array}
\right.
\end{equation}
 Notice that if $T$ is a $\mu$-block-constant tensor, then so is
$T^{(l)}_{\mbox{\rm \scriptsize out}}$ for each $l=1,...,k$. The
next theorem is Corollary~5.8 from \cite{Sendov2003b}.

\begin{theorem}
\label{dec14bc} Let $\{M_m \}$ be a sequence of symmetric matrices
converging to 0, such that $M_m/\|M_m \|$ converges to $M$. Let
$\mu$ be in $\R^n_{\downarrow}$ and $U_m \rightarrow U \in O^n$ be
a sequence of orthogonal matrices such that
\begin{equation*}
\mbox{\rm Diag}\, \mu + M_m = U_m \big( \mbox{\rm Diag}\,
\lambda(\mbox{\rm Diag}\, \mu + M_m) \big) U_m^T, \, \, \mbox{ for
all } \, \, m=1,2,....
\end{equation*}
Then for every block-constant $k$-tensor $T$ on $\R^n$, and any
permutation $\sigma$ on $\N_k$ we have
\begin{equation}
\label{kr} \lim_{m \rightarrow \infty} \frac{U_m (\Diag^{\sigma}
T) U_m^T - \Diag^{\sigma} T }{\|M_m\|} = \sum_{l=1}^k
\big(\Diag^{\sigma_{(l)}} T^{(l)}_{\mbox{\rm \scriptsize
out}}\big)[M].
\end{equation}
\end{theorem}

Again, for a fixed vector $\mu \in \R^n$, we define $k$ linear
maps
$$
T \in T^{k,n} \rightarrow T^{(l)}_{\mbox{\rm \scriptsize in}} \in
T^{k+1,n}, \mbox{ for } l=1,2,...,k,
$$
as follows:
\begin{equation}
\label{dec14c:eqn} \big(T^{\tau_l}_{\mbox{\rm \scriptsize
in}}\big)^{i_1...i_ki_{k+1}} = \left\{
\begin{array}{ll}
T^{i_1...i_{l-1}i_{k+1}i_{l+1}...i_k}, & \mbox{ if } i_l \sim
i_{k+1}, \\
0, & \mbox{ if } i_l \not\sim i_{k+1}.
\end{array}
\right.
\end{equation}
Notice that if $T$ is a block-constant tensor, then so is
$T^{\tau_l}_{\mbox{\rm \scriptsize in}}$ for each $l=1,...,k$.
Finally, we define
\begin{equation}
\label{perm-lift} \big(T^{\tau_l}\big)^{i_1...i_ki_{k+1}} =
\left\{
\begin{array}{ll}
T^{i_1...i_{l-1}i_{l}i_{l+1}...i_k}, & \mbox{ if } i_l =
i_{k+1}, \\
0, & \mbox{ if } i_l \not= i_{k+1}.
\end{array}
\right.
\end{equation}
In other words, $T^{\tau_l}$ is a $(k+1)$-tensor with entries off
the hyper plane $i_l=i_{k+1}$ equal to zero. On the hyper plane
$i_l=i_{k+1}$ we have placed the original tensor $T$. The next
theorem is Corollary~5.6 from \cite{Sendov2003b}.

\begin{theorem}
\label{jan11a} Let $U \in O^n$ be a block-diagonal orthogonal
matrix and let $\sigma$ be a permutation on $\N_k$. Let $M$ be an
arbitrary symmetric matrix, and $h \in \R^n$ be a vector, such
that $U^TM_{\mbox{\rm \scriptsize in}}U = \Diag h$. Then
\begin{enumerate}
\item for any block-constant $(k+1)$-tensor $T$ on $\R^n$,
$$
U\big(\Diag^{\sigma} (T[h]) \big)U^T =
\big(\Diag^{\sigma_{(k+1)}}T \big) [M];
$$
\item for any block-constant $k$-tensor $T$ on $\R^n$
$$
U \big(\Diag^{\sigma} (T^{\tau_l}[h]) \big)U^T =
\big(\Diag^{\sigma_{(l)}} T^{\tau_l}_{\mbox{\rm \scriptsize in}}
\big) [M], \,\,\, \mbox{ for all } l=1,...,k,
$$
\end{enumerate}
where the permutations $\sigma_{(l)}$, for $l \in \N_k$, are
defined by (\ref{permdefn}).
\end{theorem}

\section{Several standing assumptions}
\label{stndassum}

Our approach is to successively differentiate the composition $f
\circ \lambda$ where at every step we use the tensorial language
presented in Section~\ref{NABR} to simplify the calculation. More
precisely, we will define $k$-tensor valued maps
$\mathcal{A}_{\sigma} : \R^n \rightarrow T^{k,n}$, $\sigma \in
P^k$, (only in terms of the function $f$ and its partial
derivatives) such that
\begin{equation}
\label{conj:eqn} \nabla^k (f \circ \lambda)(X) = V \Big(
\sum_{\sigma \in P^k} \Diag^{\sigma}
\mathcal{A}_{\sigma}(\lambda(X)) \Big) V^T,
\end{equation}
where $X=V (\Diag \lambda(X)) V^T$. The formula for the gradient
(the case $k=1$) was originally derived in \cite{Lewis:1994b}, see
also Subsection~\ref{thegrad} below.  We showed in
\cite[Section~5]{Sendov2003a}, that having derived that formula
for $k=1$, then for $k \ge 2$ it is enough to show it under the
following three assumptions.
\begin{itemize}
\item The matrix $X$ is diagonal, $\Diag \mu$, for some vector
$\mu \in \R^n_{\downarrow}$. \\[-0.8cm]
\item The sequence $\{M_m \}$ of symmetric matrices converges to 0
and is such that $M_m/\|M_m \|$ converges to $M$. \\[-0.8cm]
\item A sequence of orthogonal matrices $U_m \in O^n$ is chosen
such that
\begin{equation*}
\mbox{\rm Diag}\, \mu + M_m = U_m \big( \mbox{\rm Diag}\,
\lambda(\mbox{\rm Diag}\, \mu + M_m) \big) U_m^T, \, \, \mbox{ for
all } \, \, m=1,2,....
\end{equation*}
and $U_m$ approaches $U \in O^n$ as $m$ goes to infinity. ($U$ is
block diagonal with blocks determined by $\mu$.)
\end{itemize}

The next lemma (the proof is a simple combination of Lemma~5.10 in
\cite{Lewis:1996} and Theorem~3.12 in \cite{HUYe:1995}) justifies
the notation that follows. Recall that $\mu \in \R^n$ partitions
$\N_n$ into $r$ blocks $I_1$,...,$I_r$.
\begin{lemma}
\label{lemma:1} For any $\mu \in \R^n_{\downarrow}$ and sequence
of symmetric matrices $M_m \rightarrow 0$ we have that
\begin{equation}
\label{dec15a:eqn} \lambda(\mbox{\rm Diag}\, \mu + M_m)^T = \mu^T
+ \left( \lambda(X_1^TM_m X_1)^T,...,\lambda(X_r^TM_mX_r)^T
\right)^T + o(\|M_m \|),
\end{equation}
where $X_l := [e^{i} \, | \, i \in I_l], \, \mbox{ for all } \,
l=1,...,r$.
\end{lemma}
Throughout the whole paper, we denote
\begin{equation}
\label{dec15b:eqn} h_m := \left( \lambda(X_1^TM_m
X_1)^T,...,\lambda(X_r^TM_mX_r)^T \right)^T.
\end{equation}
If also $M_m/\|M_m\|$ converges to $M$ as $m$ goes to infinity,
since the eigenvalues are continuous functions, we can define
\begin{equation}
\label{limh} h := \lim_{m \rightarrow \infty} \frac{h_m}{\|M_m\|}
= \left( \lambda(X_1^TM X_1)^T,...,\lambda(X_r^TMX_r)^T \right)^T.
\end{equation}
We reserve the symbols $h_m$ and $h$ to denote the above two
vectors throughout the paper. With this notation
Lemma~\ref{lemma:1} says that if $M_m \rightarrow 0$, then
\begin{equation}
\label{kamen} \lambda(\Diag \mu + M_m)^T = \mu^T + h_m +
o(\|M_m\|).
\end{equation}
If, for the fixed vector $\mu \in \R^n_{\downarrow}$, we define
\begin{align*}
M^{ij}_{\mbox{\scriptsize \rm in}} &= \left\{ \begin{array}{ll}
M^{ij}, &
\mbox{ if } i \sim j, \\
0, & \mbox{ otherwise, }
\end{array} \right.
\end{align*}
then Theorem~4.2 in \cite{Sendov2003b} says that the orthogonal
matrix $U$ is block-diagonal and satisfies
\begin{equation}
\label{diag} U^T M_{\mbox{\scriptsize \rm in}} U = \Diag h.
\end{equation}

\section{Analyticity of isolated eigenvalues}

Let $A$ be in $S^n$ and suppose that the $j$-th largest eigenvalue
is isolated, that is
$$
\lambda_{j-1}(A) > \lambda_{j}(A) > \lambda_{j+1}(A).
$$
The goal of this section is to give two justifications of the
known fact that $\lambda_{j}(\cdot)$ is an analytic function in a
neighbourhood of $A$. We call a function of several real variables
{\it analytic} at a point if in a neighbourhood of this point it
has an power series expansion. The corresponding complex variable
notion is called {\it holomorphic}.

The first justification below is from
\cite[Theorem~2.1]{TsingFanVerriest:1994}.

\begin{theorem}
\label{AnalyticCom:t} Suppose $f : \R^n \rightarrow \R$ is a
function analytic at the point $\lambda(A)$ for some $A$ in $S^n$.
Suppose also $f(Px)=f(x)$ for every permutation matrix, P, for
which $P \lambda(A) = \lambda(A)$.  Then, the function $f \circ
\lambda$ is analytic at $A$. \hfill \qed
\end{theorem}

To see how this theorem implies the analyticity of
$\lambda_{j}(\cdot)$ take
$$
f(x_1,...,x_n) = \mbox{the $j^{\mbox{\rm \scriptsize th}}$ largest
element of } \{x_1,...,x_n \}.
$$
The function $f$ is a piece-wise affine function. Moreover, for
any $x \in \R^n$ in a neighbourhood of the vector $\lambda(A)$ it
is given by
$$
f(x) = x_j.
$$
Thus, $f$ is analytic in that neighbourhood. Next, $f$ is a
symmetric function and thus by definition $f(Px)=f(x)$ for every
$x \in R^n$ and every permutation matrix $P$. Therefore by the
theorem $\lambda_j = f \circ \lambda$ is an analytic function.

For the second justification we use the following result
\cite{Arnold71}. (In the theorem below, $\lambda_i(X)$ denotes an
arbitrary eigenvalue of a matrix $X$, not necessarily the $i$'th
largest one.)

\begin{theorem}[Arnold~1971]
Suppose that $A \in \C^{n \times n}$ has $q$ eigenvalues
$\lambda_1(A),...,\lambda_q(A)$ (counting multiplicities) in an
open set $\Omega \subset \C$, and the rest $n-q$ eigenvalues not
in the closure of $\Omega$. Then, there is a neighbourhood
$\Delta$ of $A$ and holomorphic mappings $S:\Delta \rightarrow
\C^{q \times q}$ and $T:\Delta \rightarrow \C^{(n-q) \times
(n-q)}$ such that for all $X \in \Delta$
$$
X \mbox{ is similar to } \left(
\begin{array}{cc}
S(X) & 0 \\
0 & T(X)
\end{array}
\right),
$$
and $S(A)$ has eigenvalues $\lambda_1(A),...,\lambda_q(A)$. \hfill
\qed
\end{theorem}
To deduce the result we need, since the $j^{\mbox{\rm \scriptsize
th}}$ largest eigenvalue is isolated, we can find an open set
$\Omega \subset \C$, such that only that eigenvalue is in $\Omega$
and the remaining $n-1$ are not in the closure of $\Omega$. By the
theorem, there is a neighbourhood $\Delta$ of $A$ and holomorphic
mapping $S:\Delta \rightarrow \C$ such that $S(X)$ is equal to the
$j^{\mbox{\rm \scriptsize th}}$ largest eigenvalue of $X$ for all
$X$ in $\Delta$.

If $A$ is a real symmetric matrix, then the intersection of
$\Delta$ with $S^n$ is a neighbourhood of $A$ in $S^n$. Let
$\tilde{S}(X)$ denote the restriction of $S(X)$ to $\Delta \cap
S^n$. Clearly, $\tilde{S}(X)$ is holomorphic, real valued
function. Therefore, (it is a standard result in complex analysis)
the coefficients in the power series expansion of $\tilde{S}(X)$
must be real numbers. Thus, the $j^{\mbox{\rm \scriptsize th}}$
largest eigenvalue is a real analytic function in the
neighbourhood $\Delta \cap S^n$ or $A$.

All these considerations make the following observation clear.

\begin{theorem}
Suppose that $A \in S^n$ has distinct eigenvalues and $f:\R^n
\rightarrow \R$ is $k$-times (continuously) differentiable in a
neighbourhood of $\lambda(A)$. Then, $f \circ \lambda$ is
$k$-times (continuously) differentiable in a neighbourhood of $A$.
\end{theorem}

\section{The $k^{\mbox{\rm \scriptsize th}}$ derivative of functions of eigenvalues
at a matrix with distinct eigenvalues}
\label{k-thderiv}

Let $f:\R^n \rightarrow \R$ be an arbitrary $k$-times
(continuously) differentiable function. In this section, we do not
assume that $f$ is a symmetric function. Our goal in this section
is to derive a formula for the $k^{\mbox{\rm \scriptsize th}}$
derivative of $f \circ \lambda$ on the set $\lambda^{-1}(\Omega)$,
where
\begin{align*}
\Omega &= \{x \in \R^n \, | \, x_i \not= x_j \mbox{ for every } i
\not= j\}, \mbox{ and }\\
\lambda^{-1}(\Omega) &= \{A \in S^n \, | \, \lambda(A) \in \Omega
\}.
\end{align*}
Clearly $\Omega$ is a dense, open subset of $\R^n$ and
$\lambda^{-1}(\Omega)$ is a dense, open subset of $S^n$.

As an example of how one can differentiate $f \circ \lambda$, let
us consider the general situation. Let $X$, $Y$, and $Z$ be Banach
spaces and let $g:X \rightarrow Y$, $G:Y \rightarrow Z$. Then, by
applying the chain rule we have the following formulae for the
first three derivatives of $\phi = G \circ g$, (see
\cite[Section~X.4]{Bhatia:1997}) for any vectors $h_1,h_2,h_3$
from $X$:
\begin{align*}
\nabla \phi(x)[h_1] &= \nabla G(g(x))[\nabla g(x)[h_1]], \\
\nabla^2 \phi (x)[h_1, h_2] &=  \nabla^2 G(g(x))[\nabla g(x)[h_1],
\nabla g(x)[h_2]] + \nabla G(g(x))[\nabla^2 g(x)[h_1, h_2]], \\
\nabla^3 \phi(x)[h_1,h_2,h_3] &=  \nabla^3 G(g(x))[\nabla
g(x)[h_1], \nabla g(x)[h_2], \nabla g(x)[h_3]] \\
& \hspace{0.5cm} +\nabla^2 G(g(x))[\nabla g(x)[h_1], \nabla^2 g(x)[h_2, h_3]] \\
& \hspace{0.5cm} +\nabla^2 G(g(x))[\nabla g(x)[h_2], \nabla^2 g(x)[h_1, h_3]] \\
& \hspace{0.5cm} +\nabla^2 G(g(x))[\nabla g(x)[h_3], \nabla^2 g(x)[h_1, h_2]] \\
& \hspace{1cm} +\nabla G(g(x))[\nabla^3 g(x)[h_1, h_2, h_3]].
\end{align*}

In our case, we have $X= S^n$, $Y = \R^n$, $Z = \R$, $g =
\lambda$, and $G=f$. As can be seen from the above example, this
approach very quickly becomes unmanageable. The formula for the
$k$-derivative of the composition requires formulae for every
derivative of $\lambda$ up to the $k^{\mbox{\rm \scriptsize th}}$.
It is not clear how one can organize and simplify the resulting
expression into a compact, ordered formula.

Fix a vector $\mu \in \R^n_{\downarrow} \cap \Omega$. Since $\mu$
has distinct entries, every block in the partition that it defines
will have exactly one element. This means that for any $j,i \in
\N_n$, $i \sim j \Leftrightarrow i=j$, and that makes any tensor
block-constant. In particular for the matrices $X_l$, defined in
Lemma~\ref{lemma:1}, we have $X_l = [e^l], \,\,\,\, l=1,...,n$.
This implies that $h_m = \diag M_m$ and that $h=\diag M$.  Notice,
finally, how the definition of $T^{(l)}_{\mbox{\rm \scriptsize
out}}$ changes:
\begin{equation}
\label{defn-out-new} \big(T^{(l)}_{\mbox{\rm \scriptsize
out}}\big)^{i_1...i_k i_{k+1}} = \left\{
\begin{array}{ll}
0, & \mbox{ if } i_l = i_{k+1}, \\
\frac{T^{i_1...i_{l-1}i_{k+1}i_{l+1}...i_k} -
T^{i_1...i_{l-1}i_li_{l+1}...i_k}}{\mu_{i_{k+1}} - \mu_{i_l}}, &
\mbox{ if } i_l \not= i_{k+1}.
\end{array}
\right.
\end{equation}
We will derive Formula~(\ref{conj:eqn}) by induction on the order
of the derivative. For completeness, we begin by recalculating the
formula for the gradient.

\subsection{The gradient}
\label{thegrad}

Using Formulae~(\ref{kamen}) we compute
\begin{align*}
\lim_{m \rightarrow \infty} \frac{(f \circ \lambda)(\Diag \mu +
M_m) - (f \circ \lambda)(\Diag \mu)}{\|M_m\|} &= \lim_{m
\rightarrow \infty} \frac{f( \mu + h_m +
o(\|M_m\|)) - f(\mu)}{\|M_m\|} \\
&= \lim_{m \rightarrow \infty} \frac{f(\mu) + \nabla f(\mu)[h_m] +
o(\|M_m\|) - f(\mu)}{\|M_m\|} \\
&= \nabla f(\mu)[h] \\
&= \langle \nabla f(\mu), \diag M \rangle \\
&= \big(\Diag \nabla f(\mu)\big)[M].
\end{align*}
This shows that $\nabla (f \circ \lambda)(\Diag \mu) =
\Diag^{\mbox{\tiny $(1)$}} \nabla f(\mu)$. It is easy to see now
that
\begin{equation}
\label{dec14:eqn} \nabla (f \circ \lambda)(X) =
V\big(\Diag^{\mbox{\tiny $(1)$}} \nabla f(\lambda(X))\big)V^T =
V\Big(\sum_{\sigma \in P^1} \Diag^{\sigma}
\mathcal{A}_{\sigma}(\lambda(X))\Big)V^T,
\end{equation}
where $X=V(\Diag \lambda(X))V^T$ and $\mathcal{A}_{(1)}(x) =
\nabla f(x)$. Trivially, if $f$ is $k$-times (continuously)
differentiable, then $\mathcal{A}_{(1)}(x) = \nabla f(x)$ is
$(k-1)$-times (continuously) differentiable.

Note that when the eigenvalues of $X$ are not distinct, the
calculation of the gradient of $f \circ \lambda$ is almost
identical and leads to the same final formula. Indeed, using
Equation~(\ref{diag}),
$$
\nabla f(\mu)[h]= \langle \nabla f(\mu), \diag
(U^TM_{\mbox{\scriptsize \rm in}}U) \rangle = \big(U(\Diag \nabla
f(\mu))U^T \big)[M] = \big(\Diag \nabla f(\mu) \big)[M],
$$
where in the last equality we used the fact the $U$ is
block-diagonal, orthogonal and $\nabla f(\mu)$ is block-constant.

\subsection{The induction step}

Suppose now that for some $1 \le s < k$
$$
\nabla^s (f \circ \lambda)(X) = V\big(\sum_{\sigma \in P^s}
\Diag^{\sigma} \mathcal{A}_{\sigma}(\lambda(X)) \big)V^T,
$$
where $X=V(\Diag \lambda(X))V^T$. Suppose also that for every
$\sigma \in P^s$, the $s$-tensor valued map $\mathcal{A}_{\sigma}
: \R^n \rightarrow T^{s,n}$, is $(k-s)$-times (continuously)
differentiable.

Using Formulae~(\ref{kamen}), we differentiate $\nabla^s (f \circ
\lambda)$ at the matrix $\Diag \mu$:
\begin{align*}
\nabla^{s+1} (f \circ \lambda)&(\Diag \mu)[M] \\
&=\lim_{m \rightarrow \infty} \frac{\nabla^s(f \circ
\lambda)(\Diag
\mu+M_m) -\nabla^s(f \circ \lambda)(\Diag \mu)}{\|M_m\|} \\[0.1cm]
&= \lim_{m \rightarrow \infty} \frac{U_m \big( \sum_{\sigma \in
P^s} \Diag^{\sigma} \mathcal{A}_{\sigma}(\lambda(\Diag \mu+M_m))
\big)U_m^T -\sum_{\sigma \in
P^s} \Diag^{\sigma} \mathcal{A}_{\sigma}(\mu)}{\|M_m\|} \\[0.1cm]
&=\lim_{m \rightarrow \infty} \frac{\sum_{\sigma \in P^s} \Big(
U_m \big(\Diag^{\sigma} \mathcal{A}_{\sigma}(\lambda(\Diag
\mu+M_m)) \big)U_m^T -\Diag^{\sigma} \mathcal{A}_{\sigma}(\mu) \Big)}{\|M_m\|} \\[0.1cm]
&=\lim_{m \rightarrow \infty} \frac{\sum_{\sigma \in P^s} \Big(
U_m \big(\Diag^{\sigma} \mathcal{A}_{\sigma}( \mu+h_m +
o(\|M_m\|)) \big)U_m^T
-\Diag^{\sigma} \mathcal{A}_{\sigma}(\mu) \Big)}{\|M_m\|} \\[0.1cm]
&=\lim_{m \rightarrow \infty} \frac{\sum_{\sigma \in P^s} \Big(
U_m \big(\Diag^{\sigma} \big( \mathcal{A}_{\sigma}(\mu) + \nabla
\mathcal{A}_{\sigma}(\mu)[h_m] + o(\|M_m\|) \big) \big)U_m^T
-\Diag^{\sigma} \mathcal{A}_{\sigma}(\mu) \Big)}{\|M_m\|} \\[0.1cm]
&=\lim_{m \rightarrow \infty} \sum_{\sigma \in P^s} \frac{U_m
\big(\Diag^{\sigma} \mathcal{A}_{\sigma}(\mu) \big)U_m^T
-\Diag^{\sigma} \mathcal{A}_{\sigma}(\mu)}{\|M_m\|} + \sum_{\sigma
\in P^s} U \big( \Diag^{\sigma} \big( \nabla
\mathcal{A}_{\sigma}(\mu)[h] \big) \big) U^T.
\end{align*}
By Theorem~\ref{dec14bc}, since for every $\sigma \in P^s$, the
tensor $\mathcal{A}_{\sigma}(\mu)$ is block-constant, we have
\begin{align*}
\lim_{m \rightarrow \infty} \frac{U_m \big(\Diag^{\sigma}
\mathcal{A}_{\sigma}(\mu) \big)U_m^T -\Diag^{\sigma}
\mathcal{A}_{\sigma}(\mu)}{\|M_m\|} &= \sum_{l=1}^s \big(
\Diag^{\sigma_{(l)}} (\mathcal{A}_{\sigma}(\mu))^{(l)}_{\mbox{\rm
\scriptsize out}} \big) [M].
\end{align*}
By Theorem~\ref{jan11a}, since for every $\sigma \in P^s$ $\nabla
\mathcal{A}_{\sigma}(\mu)$ is a block-constant $(s+1)$-tensor, we
have
\begin{align*}
U \big(\Diag^{\sigma} \big( \nabla \mathcal{A}_{\sigma}(\mu)[h]
\big) \big) U^T = \big(\Diag^{\sigma_{(s+1)}} \nabla
\mathcal{A}_{\sigma}(\mu) \big) [M].
\end{align*}
Thus we define
\begin{align*}
\mathcal{A}_{_{\sigma_{(l)}}} &:=
(\mathcal{A}_{\sigma}(\mu))^{(l)}_{\mbox{\rm \scriptsize out}},
\, \mbox{ for all } l \in \N_s, \mbox{ and } \\
\mathcal{A}_{_{\sigma_{(s+1)}}} &:= \nabla
\mathcal{A}_{\sigma}(\mu).
\end{align*}
Putting everything together and conclude that for every symmetric
matrix $M$:
$$
\nabla^{s+1} (f \circ \lambda)(\Diag \mu)[M] = \Big(
\hspace{-0.3cm} \sum_{\substack{ \sigma \in P^s \\
\,\,\,\,\,\,\,\, l \in \N_{s+1}}} \hspace{-0.3cm}
\Diag^{\sigma_{(l)}} \mathcal{A}_{\sigma_{(l)}}(\mu) \Big)[M].
$$
Notice the parameters of the summation sign in the above formula.
As $\sigma$ goes over the elements of $P^s$ and as $l$ goes over
the set $\N_{s+1}$ the permutation $\sigma_{(l)}$ covers, in a
one-to-one manner, all permutations in $P^{s+1}$.  Now, the
comments in \cite[Section~5]{Sendov2003a} show that
$$
\nabla^{s+1} (f \circ \lambda)(X) = V \Big( \hspace{-0.3cm}
\sum_{\substack{ \sigma \in P^s \\ \,\,\,\,\,\,\,\, l \in
\N_{s+1}}} \hspace{-0.3cm} \Diag^{\sigma_{(l)}}
\mathcal{A}_{\sigma_{(l)}}(\lambda(X))\Big) V^T,
$$
where $X=V(\Diag \lambda(X))V^T$.

To finish the induction, we have to show that the $(s+1)$-tensor
valued maps $\mathcal{A}_{\sigma_{(l)}}(\cdot)$ are at least
$(k-s-1)$-times (continuously) differentiable. This is clear when
$l=s+1$ and $\sigma \in P^s$, since $\mathcal{A}_{\sigma}(\cdot)$
is $(k-s)$-times (continuously) differentiable for every $\sigma
 \in P^s$. For the rest of the maps this is also easy to see. Every entry
in $\mathcal{A}_{\sigma_{(l)}}$ is the difference of two entries
of $\mathcal{A}_{\sigma}$ divided by a quantity that never becomes
zero over the set $\Omega$. This shows that over the set $\Omega$,
$\mathcal{A}_{\sigma_{(l)}}(\cdot)$ is $(k-s)$-times
(continuously) differentiable for every $\sigma \in P^s$ and every
$l \in \N_s$.

We summarize everything in the next theorem.

\begin{theorem}
\label{dist-main}
 Let $X$ be a symmetric matrix with distinct
eigenvalues and let $f$ be a $k$-times (continuously)
differentiable function on $\R^n$. Let $V$ be an orthogonal matrix
such that $X=V(\Diag \lambda(X))V^T$. Then, $f \circ \lambda$ is
$k$-times (continuously) differentiable function at $X$. Moreover
if $\nabla^{s}(f \circ \lambda)$, for some $s < k$, is given by
$$
\nabla^s (f \circ \lambda)(X) = V\Big( \sum_{\sigma \in P^{s}}
\Diag^{\sigma} \mathcal{A}_{\sigma}(\lambda(X)) \Big)V^T,
$$
for some $s$-tensor valued mappings $\mathcal{A}_{\sigma} : \R^n
\rightarrow T^{s,n}$, for every $\sigma \in P^{s}$, then
$\nabla^{(s+1)}(f \circ \lambda)$ is given by
\begin{equation}
\label{dist-main-form} \nabla^{(s+1)} (f \circ \lambda)(X) =
V\Big(\hspace{-0.3cm} \sum_{\substack{\sigma \in P^{s} \\
\,\,\,\,\,\,\,\, l \in \N_{s+1}}} \hspace{-0.3cm}
\Diag^{\sigma_{(l)}} \mathcal{A}_{\sigma_{(l)}}(\lambda(X))
\Big)V^T,
\end{equation}
 where
\begin{align*}
\mathcal{A}_{\sigma_{(l)}} &= \big( \mathcal{A}_{\sigma}
\big)^{(l)}_{\mbox{\rm \scriptsize out}},
\, \mbox{ for all } \, l \in \N_s, \mbox{ and } \\
\mathcal{A}_{\sigma_{(s+1)}} &= \nabla \mathcal{A}_{\sigma}.
\end{align*}
\end{theorem}

\section{The $k^{\mbox{\rm \scriptsize th}}$ derivative of
separable spectral functions}

In this section we show that Formula~(\ref{conj:eqn}) holds at an
arbitrary symmetric matrix $X$ (not necessarily with distinct
eigenvalues) for a subclass of spectral functions that we now
describe.

Let $g$ be a real function on an interval $I$. If $D=
\Diag(\lambda_1,...,\lambda_n)$ is a diagonal matrix with diagonal
entries $\lambda_i \in I$, $i=1,...,n$, we define
\begin{equation}
\label{G-defn-1} G(D) = \Diag (g(\lambda_1),...,g(\lambda_1)).
\end{equation}
If $X$ is a symmetric matrix with eigenvalues $\lambda_i$ in $I$,
we choose an orthogonal matrix $V$ such that $X = V(\Diag
\lambda(X)) V^T$ and, then define
\begin{equation}
\label{G-defn-2} G(X) = V G(\Diag \lambda(X)) V^T.
\end{equation}
In this way we obtain a (well-defined) symmetric-matrix valued
function with domain the set of all matrices $X$ with eigenvalues
in $I$.

These functions have been the object of recent interest in
optimization \cite{SunSun:2002a} and the main object of
\cite[Chapter V]{Bhatia:1997}, where their gradient is computed
using an approximation argument. Notice that $G(X)$ is just the
gradient (see Formula~(\ref{dec14:eqn})) of the spectral function
$f \circ \lambda$, where $f(x)=\tilde{g}(x_1)+ \cdots +
\tilde{g}(x_n)$, and $\tilde{g}(s) = \int_{0}^s g(t) \, dt$. That
is why we will call those functions {\it separable spectral
functions}.

\subsection{Description of the $k^{\mbox{\rm \scriptsize th}}$ derivative}

Let $g :I \rightarrow \R$ be $k$-times (continuously)
differentiable. Define the symmetric function $g^{[(12)]}(x,y) : I
\times I \rightarrow \R$ as
\begin{equation}
\label{g-(12)} g^{[(12)]}(x,y) = \left\{
\begin{array}{ll}
\displaystyle\frac{g(x) - g(y)}{x - y}, & \mbox{
if } x \not= y, \\[0.3cm]
g'(x), & \mbox{ if } x=y.
\end{array}
\right.
\end{equation}
The integral representation $g^{[(12)]}(x,y) = \int_{0}^1 g'(y +
t(x-y)) \, dt$ shows that $g^{[(12)]}(x,y)$ is as smooth, in both
arguments, as $g'$.

Denote by $\tilde{P}^k$ the set of all permutations from $P^k$
that have one cycle in their cycle decomposition. Clearly
$|\tilde{P}^k|=(k-1)!$. Notice that for every $\sigma \in
\tilde{P}^k$ and every $l \in \N_k$ we have $\sigma_{(l)} \in
\tilde{P}^{k+1}$. Moreover, as $\sigma$ varies over $\tilde{P}^k$
and $l$ varies over $\N_k$, the permutation $\sigma_{(l)}$ varies
over $\tilde{P}^{k+1}$ in a one-to-one and onto fashion.

Suppose that for every $\sigma \in \tilde{P}^k$ we have defined
the function $g^{[\sigma]}(x_1,...,x_k)$ on the set $I \times I
\times \cdots \times I$, $k$-times, and suppose that these
functions are as smooth as $g^{(k-1)}$ (the $(k-1)$-th derivative
of $g$). For every $\sigma \in \tilde{P}^k$ and every $l \in \N_k$
we define the function $g^{[\sigma_{(l)}]}(x_1,...,x_k, x_{k+1})$
as follows:
\begin{equation}
\label{ind-fn-defn} g^{[\sigma_{(l)}]}(x_1,...,x_{k+1}) = \left\{
\begin{array}{ll}
\nabla_{l} g^{[\sigma]}(x_1,...,x_k), & \mbox{ if } x_l = x_{k+1}
\\[0.2cm]
\displaystyle\frac{g^{[\sigma]}(x_1,...,x_{l},...,x_k)
-g^{[\sigma]}(x_1,...,x_{k+1},...,x_k)}{x_l-x_{k+1}}, & \mbox{ if
} x_l \not= x_{k+1},
\end{array}
\right.
\end{equation}
where in the second case of the definition, both $x_l$ and
$x_{k+1}$ are in $l$-th position, and $\nabla_l$ denotes the
partial derivative with respect to the $l$-th argument. Using the
integral formula
$$
g^{[\sigma_{(l)}]}(x_1,...,x_{k+1}) = \int_{0}^1
\nabla_lg^{[\sigma]}(x_1,...,x_{l-1},x_{k+1}+t(x_{l}-x_{k+1}),x_{l+1},...,x_k)
\, dt,
$$
for every $l \in \N_k$, we see that
$g^{[\sigma_{(l)}]}(x_1,...,x_{k+1})$ is as smooth as $g^{(k)}$,
the $k$-th derivative of $g$.

Finally, for every $s \in \{2,3,...,k+1\}$ and every $\sigma \in
\tilde{P}^s$, we define a $s$-tensor valued map
\begin{align}
\mathbf{g}^{[\sigma]}:\R^n & \,\rightarrow T^{s,n}, \mbox{ where } \nonumber \\[-0.3cm]
\label{ILJen}  \\[-0.3cm]
\big(\mathbf{g}^{[\sigma]}(\mu)\big)^{i_1...i_s} &:=
g^{[\sigma]}(\mu_{i_1},...,\mu_{i_s}). \nonumber
\end{align}
Clearly, if $(i_1,...,i_s) \sim_{\mu} (j_1,...,j_s)$, then
$\big(\mathbf{g}^{[\sigma]}(\mu)\big)^{i_1...i_s}=\big(\mathbf{g}^{[\sigma]}(\mu)\big)^{j_1...j_s}$,
which shows that $\mathbf{g}^{[\sigma]}(\mu)$ is a
$\mu$-block-constant tensor for every $\mu$. Moreover, the map
$\mathbf{g}^{[\sigma]}:\R^n \rightarrow T^{s,n}$ is still as
smooth as $g^{(s-1)}$, for every $s=2,3,...,k+1$.

We are now ready to formulate the second main result of this work.
The proof is given in the next subsection. A comparison between
Theorem~\ref{sep-main} and Theorem~\ref{dist-main} is given at the
end of Subsection~\ref{ind-step-sep}.

\begin{theorem}
\label{sep-main} Let $g$ be a $k$-times (continuously)
differentiable function defined on an interval $I$. Let $X$ be a
symmetric matrix with eigenvalues in the interval $I$, and let $V$
be an orthogonal matrix such that $X=V(\Diag \lambda(X))V^T$.
Then, the matrix valued function $G$ defined by (\ref{G-defn-1})
and (\ref{G-defn-2}) is $k$-times (continuously) differentiable at
$X$. Moreover its $k$-th derivative, $\nabla^{k} G(X)$, is given
by the formula
\begin{equation}
\label{sep-main-form} \nabla^{k} G(X) = V\Big(\hspace{0cm}
\sum_{\substack{\sigma \in \tilde{P}^{k+1}}} \hspace{0cm}
\Diag^{\sigma} \mathbf{g}^{[\sigma]}(\lambda(X)) \Big)V^T,
\end{equation}
where the $(k+1)$-tensor valued maps
$\mathbf{g}^{[\sigma]}(\cdot)$ are defined by
Equation~(\ref{ILJen}).
\end{theorem}

\subsection{Proof of Theorem~\ref{sep-main}: the gradient}

Let $X$ be an $n \times n$ symmetric matrix with all eigenvalues
in $I$ and such that $X=V(\Diag \lambda(X))V^T$ for some
orthogonal matrix $V$. The formula for the gradient of separable
spectral functions has been known for a while. For example, using
approximation techniques, it was shown in \cite{Bhatia:1997} that
for any two symmetric matrices $H_1$ and $H_2$
\begin{equation}
\label{bhatia-grad} \nabla G(X)[H_1,H_2] = \langle V
\big(g^{[(12)]}(\lambda(X)) \circ (V^TH_1V) \big)V^T, H_2 \rangle,
\end{equation}
where `$\circ$' stands for the usual Hadamard product.

In this subsection, we will give a direct derivation of the
gradient and as a result a slightly different representation of
the above formula.

For convenience we denote $g(x):=(g(x_1),...,g(x_n))^T$, and
$\nabla g(x):=(g'(x_1),...,g'(x_n))^T$, for any $x \in \R^n$. Thus
we compute:
\begin{align*}
\nabla G(\Diag \mu)[M]  &= \lim_{m \rightarrow \infty}
\frac{G(\Diag \mu+M_m)-G(\Diag \mu)}{\|M_m\|} \\
&= \lim_{m \rightarrow \infty} \frac{U_m \Diag
g(\lambda(\Diag \mu+M_m)) U_m^T -\Diag g(\mu)}{\|M_m\|} \\
&=\lim_{m \rightarrow \infty} \frac{U_m (\Diag
g(\mu+h_m+ o(\|M_m\|))) U_m^T -\Diag g(\mu)}{\|M_m\|} \\
&= \lim_{m \rightarrow \infty} \frac{U_m \big( \Diag
g(\mu) + (\Diag \nabla g(\mu))[ h_m ]+o(\|M_m\|)) \big) U_m^T -\Diag g(\mu)}{\|M_m\|} \\
&=\lim_{m \rightarrow \infty} \frac{U_m \big( \Diag g(\mu) \big)
U_m^T-\Diag g(\mu)}{\|M_m\|} + U\big((\Diag \nabla
g(\mu))[h]\big)U^T.
\end{align*}
It is important to notice that both vectors $g(\mu)$ and $\nabla
g(\mu)$ are block-constant (with respect to $\mu$).  We use the
second part of Corollary~\ref{jan11a} with $\sigma=(1)$, $k=l=1$,
(notice that $\sigma_{(l)}=(12)$) and $T=\nabla g(\mu)$ to develop
the second term above:
\begin{align*}
U\big((\Diag \nabla g(\mu))[h]\big)U^T = \big(\Diag^{\mbox{\tiny
(12)}} (\nabla g(\mu))_{\mbox{\rm \scriptsize in}}^{\mbox{\tiny
(12)}} \big)[M].
\end{align*}
We now use Corollary~\ref{dec14bc} with $k=1$, $\sigma = (1)$, and
$T=g(\mu)$ to find the limit:
\begin{align*}
\lim_{m \rightarrow \infty} \frac{U_m \big( \Diag g(\mu) \big)
U_m^T-\Diag g(\mu)}{\|M_m\|} = \big(\Diag^{\mbox{\tiny (12)}}
(g(\mu))_{\mbox{\rm \scriptsize out}}^{\mbox{\tiny (12)}} \big)
[M].
\end{align*}
Putting everything together we get
\begin{align*}
\nabla G(\Diag \mu)  &= \Diag^{\mbox{\scriptsize (12)}} (\nabla
g(\mu))_{\mbox{\rm \scriptsize in}}^{\mbox{\scriptsize (12)}}  +
 \Diag^{\mbox{\scriptsize (12)}} (g(\mu))_{\mbox{\rm \scriptsize
out}}^{\mbox{\scriptsize (12)}}  = \Diag^{\mbox{\scriptsize (12)}}
\mathbf{g}^{[(12)]}(\mu),
\end{align*}
where we used the easy to check fact that
$\mathbf{g}^{[(12)]}(\mu) = (\nabla g(\mu))_{\mbox{\rm \scriptsize
in}}^{\mbox{\scriptsize (12)}} + (g(\mu))_{\mbox{\rm \scriptsize
out}}^{\mbox{\scriptsize (12)}}$.  Now, using
Lemma~\ref{simplelem} it is easy to see the following result (when
$X$ is arbitrary symmetric matrix, not just diagonal).

\begin{theorem}
Let $g \in C^1(I)$ and let $X$ be a symmetric matrix with all
eigenvalues in $I$. Then,
\begin{equation}
\label{grad-form} \nabla G(X) = V\big(\Diag^{\mbox{\scriptsize \rm
(12)}} \mathbf{g}^{[(12)]}(\lambda(X)) \big)V^T,
\end{equation}
 where $X=V(\Diag \lambda(X))V^T$.
\end{theorem}

For the sake of completeness, we show that
Formula~(\ref{grad-form}) is indeed the same as
Formula~(\ref{bhatia-grad}). This is achieved when in the next
result one substitutes the matrix $A$ with the matrix
$\mathbf{g}^{[(12)]}(\lambda(X))$.

\begin{proposition} For any $n \times n$ matrix $A$, any orthogonal $V$, and any
symmetric $H_1$ and $H_2$, we have the equality
$$
\big(V(\Diag^{\mbox{\scriptsize \rm (12)}} A)V^T\big) [H_1, H_2] =
\langle V \big(A \circ (V^TH_1V) \big)V^T, H_2 \rangle,
$$
where `$\circ$' stands for the ordinary Hadamard product.
\end{proposition}

\begin{proof}
We develop the two sides of the stated equality and compare the
results. By Theorem~\ref{jen3a}, the left-hand side is equal to
\begin{align*}
V\big(\Diag^{\mbox{\scriptsize (12)}} A \big)V^T [H_1, H_2] &=
\langle A, \tilde{H}_1 \circ_{\hspace{-0.05cm}_{\mbox{\tiny
$(12)$}}} \tilde{H}_2 \rangle.
\end{align*}
On the other hand
\begin{align*}
\langle V \big(A \circ (V^TH_1V) \big)V^T, H_2 \rangle &= \langle
A \circ \tilde{H}_1, \tilde{H}_2 \rangle = \langle A, \tilde{H}_1
\circ \tilde{H}_2 \rangle.
\end{align*}
Finally it is easy to check directly from the definitions that
$\tilde{H}_1 \circ_{\hspace{-0.05cm}_{\mbox{\tiny $(12)$}}}
\tilde{H}_2 = \tilde{H}_1 \circ \tilde{H}_2^T = \tilde{H}_1 \circ
\tilde{H}_2$, where in the last equality we used that
$\tilde{H}_2$ is symmetric. \hfill \qed
\end{proof}

\subsection{Proof of Theorem~\ref{sep-main}: the induction step}
\label{ind-step-sep}

Suppose that $g:I \rightarrow \R$ is $k$-times (continuously)
differentiable, and that the formula for the $(s-1)$-th derivative
($2 \le s < k+1$) of $G$ at the matrix $X$ is given by
$$
\nabla^{(s-1)} G(X) = V\Big(\sum_{\sigma \in \tilde{P}^{s}}
\Diag^{\sigma} \mathbf{g}^{[\sigma]}(\lambda(X)) \Big)V^T.
$$
The $s$-tensor-valued maps $\mathbf{g}^{[\sigma]} : \R \rightarrow
T^{s,n}$ are at least $(k-s+1)$-times (continuously)
differentiable. As we explained in Section~\ref{stndassum}, it is
enough to derive the formula for $\nabla^{s} G(X)$ only in the
case when $X=\Diag \mu$ for some $\mu \in \R^n_{\downarrow}$. We
compute:
\begin{align*}
&\nabla^{s} G(\Diag \mu)[M] = \lim_{m \rightarrow \infty}
\frac{\nabla^{(s-1)} G(\Diag \mu+M_m) - \nabla^{(s-1)} G(\Diag
\mu)}{\|M_m\|} \\
&= \lim_{m \rightarrow \infty} \frac{U_m\Big(\sum_{\sigma \in
\tilde{P}^{s}} \Diag^{\sigma} \mathbf{g}^{[\sigma]}(\lambda(\Diag
\mu+M_m)) \Big)U_m^T - \sum_{\sigma \in
\tilde{P}^{s}} \Diag^{\sigma} \mathbf{g}^{[\sigma]}(\mu)}{\|M_m\|} \\
&= \lim_{m \rightarrow \infty} \frac{U_m\Big(\sum_{\sigma \in
\tilde{P}^{s}} \Diag^{\sigma}
\mathbf{g}^{[\sigma]}(\mu+h_m+o(\|M_m\|))
\Big)U_m^T - \sum_{\sigma \in \tilde{P}^{s}} \Diag^{\sigma} \mathbf{g}^{[\sigma]}(\mu)}{\|M_m\|} \\
&= \lim_{m \rightarrow \infty} \frac{U_m\Big(\sum_{\sigma \in
\tilde{P}^{s}} \Diag^{\sigma}
\big(\mathbf{g}^{[\sigma]}(\mu)+\nabla
\mathbf{g}^{[\sigma]}(\mu)[h_m]+o(\|M_m\|) \big)
\Big)U_m^T - \sum_{\sigma \in \tilde{P}^{k}} \Diag^{\sigma} \mathbf{g}^{[\sigma]}(\mu)}{\|M_m\|} \\
&= \lim_{m \rightarrow \infty} \frac{U_m\Big(\sum_{\sigma \in
\tilde{P}^{s}} \Diag^{\sigma} \mathbf{g}^{[\sigma]}(\mu)
\Big)U_m^T - \sum_{\sigma \in \tilde{P}^{s}} \Diag^{\sigma}
\mathbf{g}^{[\sigma]}(\mu)}{\|M_m\|} + U\Big(\sum_{\sigma \in
\tilde{P}^{s}} \Diag^{\sigma} \big(\nabla
\mathbf{g}^{[\sigma]}(\mu)[h]\big) \Big)U^T.
\end{align*}
First, using Theorem~\ref{dec14bc}, we wrap up the limit in the
above formula:
\begin{align}
\label{first-eqn} \lim_{m \rightarrow \infty}
\frac{U_m\Big(\sum_{\sigma \in \tilde{P}^{s}} \Diag^{\sigma}
\mathbf{g}^{[\sigma]}(\mu) \Big)U_m^T - \sum_{\sigma \in
\tilde{P}^{s}} \Diag^{\sigma} \mathbf{g}^{[\sigma]}(\mu)}{\|M_m\|}
= \sum_{\substack{\sigma \in \tilde{P}^{s} \\ \, l \in \N_s}}
\big(
\Diag^{\sigma_{(l)}}(\mathbf{g}^{[\sigma]}(\mu))^{(l)}_{\mbox{\scriptsize
\rm out}} \big)[M].
\end{align}
Next, we focus our attention on the gradient $\nabla
\mathbf{g}^{[\sigma]}(\mu)$. Using the definition,
Equation~(\ref{ILJen}), we see that
\begin{align}
\nonumber \nabla
\big[(\mathbf{g}^{[\sigma]}(\mu))^{i_1...i_s}\big] &= \sum_{l=1}^s
\nabla_l g^{[\sigma]}(\mu_{i_1},...,\mu_{i_s})e^{i_l} \\[-0.4cm]
\label{repr1} \\[-0.3cm]
&= \nonumber \sum_{l=1}^s
g^{[\sigma_{(l)}]}(\mu_{i_1},...,\mu_{i_s},\mu_{i_l})e^{i_l},
\end{align}
where for the second equality we used
Equation~(\ref{ind-fn-defn}). This prompts us to define the
$s$-tensor-valued map
\begin{align}
T_{l}:\R^n & \, \rightarrow T^{s,n}, \mbox{ where }\nonumber \\[-0.3cm]
\label{new-tens} \\[-0.3cm]
(T_{l}(\mu))^{i_1...i_s} &:=
g^{[\sigma_{(l)}]}(\mu_{i_1},...,\mu_{i_s},\mu_{i_l}). \nonumber
\end{align}
for every $l \in \N_s$. Notice that $T_{l}(\mu)$ is a
$\mu$-block-constant $s$-tensor, for every $\mu$ and every $l \in
\N_s$.

\begin{lemma}
\label{decomp-pieces} The gradient of $\mathbf{g}^{[\sigma]}(\mu)$
allows the following decomposition
\begin{equation}
\label{decomps-1} \nabla \mathbf{g}^{[\sigma]}(\mu) = \sum_{l=1}^s
\big(T_{l}(\mu)\big)^{\tau_l},
\end{equation}
where the ``lifting'' $\big(T_{l}(\mu)\big)^{\tau_l}$ is defined
by Equation~(\ref{perm-lift}).
\end{lemma}

\begin{proof}
Fix a multi index $(i_1,...,i_s)$. By definition of the gradient
$\nabla g^{[\sigma]}(\mu)$ we have that
$$
\big((\nabla \mathbf{g}^{[\sigma]}(\mu))^{i_1...i_s,1},(\nabla
\mathbf{g}^{[\sigma]}(\mu))^{i_1...i_s,2},...,(\nabla
\mathbf{g}^{[\sigma]}(\mu))^{i_1...i_s,n}\big)^T = \nabla
\big[(\mathbf{g}^{[\sigma]}(\mu))^{i_1...i_s}\big].
$$
We compute the $p$-th entry in the above vector. On one hand,
using Equation~(\ref{repr1}), we get:
\begin{align*}
(\nabla \mathbf{g}^{[\sigma]}(\mu))^{i_1...i_s,p} &= \sum_{\substack{l=1 \\
i_l=p}}^s g^{[\sigma_{(l)}]}(\mu_{i_1},...,\mu_{i_s},\mu_{i_l}).
\end{align*}
On the other, using Equation~(\ref{perm-lift}), we evaluate the
right-hand side of (\ref{decomps-1}):
\begin{align*}
\Big(\sum_{l=1}^s \big(T_{l}(\mu)\big)^{\tau_l}\Big)^{i_1...i_s,p}
&= \sum_{l=1}^s
\big(\big(T_{l}(\mu)\big)^{\tau_l}\big)^{i_1...i_s,p} \\
&= \sum_{l=1}^s \big(T_{l}(\mu)\big)^{i_1...i_s}\delta_{i_lp} \\
&= \sum_{\substack{l=1 \\ i_l=p}}^s \big(T_{l}(\mu)\big)^{i_1...i_s} \\
&= \sum_{\substack{l=1 \\
i_l=p}}^s g^{[\sigma_{(l)}]}(\mu_{i_1},...,\mu_{i_s},\mu_{i_l}).
\\[-1.3cm]
\end{align*}
\hfill \qed
\end{proof}

We now continue the evaluation of $\nabla^s G(\Diag \mu)[M]$.
Using Theorem~\ref{jan11a} in the last equality below, we find
that
\begin{align}
U\Big(\sum_{\sigma \in \tilde{P}^{s}} \Diag^{\sigma} \big(\nabla
g^{[\sigma]}(\mu)[h]\big) \Big)U^T &= U\Big(\sum_{\sigma \in
\tilde{P}^{s}} \Diag^{\sigma} \Big(\Big(\sum_{l=1}^s
\big(T_{l}(\mu)\big)^{\tau_l} \Big)[h]\Big) \Big)U^T \nonumber \\
&= U\Big(\sum_{\sigma \in \tilde{P}^{s}} \Diag^{\sigma}
\Big(\sum_{l=1}^s \big(T_{l}(\mu)\big)^{\tau_l}[h] \Big) \Big)U^T \nonumber \\
&= U\Big(\sum_{\sigma \in \tilde{P}^{s}} \sum_{l=1}^s
\Diag^{\sigma} \big( \big(T_{l}(\mu)\big)^{\tau_l}[h] \big)
\Big)U^T \nonumber \\
&= \sum_{\sigma \in \tilde{P}^{s}} \sum_{l=1}^s U
\Big(\Diag^{\sigma} \big( \big(T_{l}(\mu)\big)^{\tau_l}[h] \big)
\Big)U^T \nonumber \\
&= \label{second-eqn} \sum_{\sigma \in \tilde{P}^{s}} \sum_{l=1}^s
\Big(\Diag^{\sigma_{(l)}}
\big(T_{l}(\mu)\big)^{\tau_l}_{\mbox{\scriptsize \rm in}}
\Big)[M].
\end{align}
This already shows that $\nabla^{(s-1)} G(\Diag \mu)$ is
differentiable. All that is left to do, now, is to show that
$\nabla^{s} G(\Diag \mu)$ has the desired form and properties. The
last step is formulated in the next lemma.

\begin{lemma}
\label{lastpush} For every $l \in \N_s$ the following identity
holds
$$
\mathbf{g}^{[\sigma_{(l)}]}(\mu) =
\big(T_{l}(\mu)\big)^{\tau_l}_{\mbox{\scriptsize \rm in}} +
\big(\mathbf{g}^{[\sigma]}(\mu)\big)^{(l)}_{\mbox{\scriptsize \rm
out}}.
$$
\end{lemma}

\begin{proof}
Fix a number $l \in \N_s$ and a multi index
$(i_1,...,i_s,i_{s+1})$. We consider two cases depending on
whether or not $\mu_{i_{s+1}}$ equals $\mu_{i_l}$.

{\bf Case I.} Suppose $i_l \sim_{\mu} i_{s+1}$. Then, the entry of
the left-hand side, corresponding to the multi index
$(i_1,...,i_s,i_{s+1})$ is
\begin{align*}
\big(\mathbf{g}^{[\sigma_{(l)}]}(\mu)\big)^{i_1...i_si_{s+1}} &=
g^{[\sigma_{(l)}]}(\mu_{i_1},...,\mu_{i_s},\mu_{i_{s+1}}) =
\nabla_l g^{[\sigma]}(\mu_{i_1},...,\mu_{i_s}).
\end{align*}
On the other hand, the right-hand side evaluates to
\begin{align*}
\big(\big(T_{l}(\mu)\big)^{\tau_l}_{\mbox{\scriptsize \rm in}} +
\big(\mathbf{g}^{[\sigma]}(\mu)\big)^{(l)}_{\mbox{\scriptsize \rm
out}}\big)^{i_1...i_si_{s+1}} &=
\big(\big(T_{l}(\mu)\big)^{\tau_l}_{\mbox{\scriptsize \rm in}}
\big)^{i_1...i_si_{s+1}} +
\big(\big(\mathbf{g}^{[\sigma]}(\mu)\big)^{(l)}_{\mbox{\scriptsize
\rm out}}\big)^{i_1...i_si_{s+1}} \\
&= \big(\big(T_{l}(\mu)\big)^{\tau_l}_{\mbox{\scriptsize \rm in}}
\big)^{i_1...i_si_{s+1}} + 0 \\
&= \big(T_{l}(\mu)\big)^{i_1...i_s} \\
&= g^{[\sigma_{(l)}]}(\mu_{i_1},...,\mu_{i_s},\mu_{i_l}) \\
&= \nabla_l g^{[\sigma]}(\mu_{i_1},...,\mu_{i_s}),
\end{align*}
where in the third equality we used Equation~(\ref{dec14c:eqn})
and the fact that $T_{l}(\mu)$ is block-constant.

{\bf Case II.} Suppose $i_l \not\sim_{\mu} i_{s+1}$. Then, the
entry of the left-hand side, corresponding to the multi index
$(i_1,...,i_s,i_{s+1})$ is
\begin{align*}
\big(\mathbf{g}^{[\sigma_{(l)}]}(\mu)\big)^{i_1...i_si_{s+1}} &=
g^{[\sigma_{(l)}]}(\mu_{i_1},...,\mu_{i_s},\mu_{i_{s+1}}) \\
&= \frac{g^{[\sigma]}(\mu_{i_1},...,\mu_{i_l},...,\mu_{i_s})
-g^{[\sigma]}(\mu_{i_1},...,\mu_{i_{s+1}},...,\mu_{i_s})}{\mu_{i_l}-\mu_{i_{s+1}}},
\end{align*}
where both $\mu_{i_l}$ and $\mu_{i_{s+1}}$ are in position $l$. On
the other hand, the right-hand side evaluates to
\begin{align*}
\big(\big(T_{l}(\mu)\big)^{\tau_l}_{\mbox{\scriptsize \rm in}} +
\big(\mathbf{g}^{[\sigma]}(\mu)\big)^{(l)}_{\mbox{\scriptsize \rm
out}}\big)^{i_1...i_si_{s+1}} &=
\big(\big(T_{l}(\mu)\big)^{\tau_l}_{\mbox{\scriptsize \rm in}}
\big)^{i_1...i_si_{s+1}} +
\big(\big(\mathbf{g}^{[\sigma]}(\mu)\big)^{(l)}_{\mbox{\scriptsize
\rm out}}\big)^{i_1...i_si_{s+1}} \\
&=  0 +
\big(\big(\mathbf{g}^{[\sigma]}(\mu)\big)^{(l)}_{\mbox{\scriptsize
\rm out}}\big)^{i_1...i_si_{s+1}}\\
&=\frac{\big(\mathbf{g}^{[\sigma]}(\mu)\big)^{i_1...i_{l-1}i_{s+1}i_{l+1}...i_s}
-\big(\mathbf{g}^{[\sigma]}(\mu)\big)^{i_1...i_{l-1}i_li_{l+1}...i_s}}{\mu_{i_{s+1}}-\mu_{i_l}} \\
&=\frac{g^{[\sigma]}(\mu_{i_1},...,\mu_{i_{s+1}},...,\mu_{i_s})
-g^{[\sigma]}(\mu_{i_1},...,\mu_{i_l},...,\mu_{i_s})}{\mu_{i_{s+1}}-\mu_{i_l}}.
\end{align*}
In both cases, the two sides are equal and we are done. \hfill
\qed
\end{proof}

Putting Equations~(\ref{first-eqn}) and (\ref{second-eqn})
together, and using Lemma~\ref{lastpush} concludes the inductive
step and proves Theorem~\ref{sep-main}.

In the special case when matrix $X$ has distinct eigenvalues, it
seems that Theorem~\ref{dist-main} and Theorem~\ref{sep-main} give
two different formulae for the higher-order derivatives of a
separable spectral function. We now reconcile the differences.
Suppose we have the formula for $\nabla^s G(X)$ given by
Equation~(\ref{sep-main-form}) and apply to it the inductive
procedure described in Theorem~\ref{dist-main} to obtain
$\nabla^{(s+1)} G(X)$. The calculations in
Subsection~\ref{ind-step-sep} showed that the gradient
$\mathcal{A}_{\sigma_{(s+1)}} = \nabla \mathcal{A}_{\sigma}$ can
be partitioned into $s$ pieces (Lemma~\ref{decomp-pieces}) and
each piece can be added as an $s$-dimensional ``diagonal plane''
(Lemma~\ref{lastpush}) to a corresponding tensor
$\mathcal{A}_{\sigma_{(l)}}$ for $l \in \N_s$. Doing that, we will
arrive at the formula for $\nabla^{(s+1)} G(X)$ given by
Theorem~\ref{sep-main}.

\subsection{$C^k$ separable spectral functions}

Theorem~\ref{sep-main} holds for every $k$-times differentiable
functions $g$. If in addition $g$ in $k$-times continuously
differentiable, then Formula~(\ref{sep-main-form}) can be
 significantly simplified. This is what we will describe in this
section.  In particular, we will show three properties of the
functions $g^{[\sigma]}(x_1,...,x_{s})$, for every $2 \le s \le
k+1$ and every $\sigma \in \tilde{P}_{s}$.  First, we will give a
compact determinant formula for computing
$g^{[\sigma]}(x_1,...,x_{s})$ directly. Second, as a consequence
of the determinant formula we will see that
$g^{[\sigma]}(x_1,...,x_{s})$ is a symmetric function on its $s$
arguments. Finally, third, denoting $\sigma_s = (12...s)$, all
functions $g^{[\sigma]}(x_1,...,x_{s})$ can be obtained from
$g^{[\sigma_s]}(x_1,...,x_{s})$ by a permutation of its arguments.
(Thus, knowing one of the tensors in
Formula~(\ref{sep-main-form}), namely
$\mathbf{g}^{[\sigma_s]}(\mu)$, we can obtain the rest by
permuting its ``rows'' and ``columns''.)

Denote by $V(x_1,...,x_s)$ the Vandermonde determinant
$$
V(x_1,...,x_s) = \left| \begin{array}{cccc} x_1^{s-1} & x_2^{s-1} & \cdots & x_s^{s-1} \\
\vdots & \vdots & \ddots & \vdots \\
x_1 & x_2 & \cdots & x_s \\
1 & 1 & \cdots & 1
\end{array} \right| = \prod_{j<i} (x_j-x_i).
$$
For any $y \in \R^s$, denote by $V\big(\substack{y_1,...,y_s \\
x_1,...,x_s}\big)$ the determinant
$$
V\big(\substack{y_1,...,y_s \\ x_1,...,x_s}\big) = \left|
\begin{array}{cccc} y_1 & y_2 & \cdots
& y_s \\
x_1^{s-2} & x_2^{s-2} & \cdots & x_s^{s-2} \\
\vdots & \vdots & \ddots & \vdots \\
x_1 & x_2 & \cdots & x_s \\
1 & 1 & \cdots & 1
\end{array} \right|.
$$

\begin{lemma}
\label{vand-ident} For any vector $(x_1,...,x_s, x_{s+1})$ with
distinct coordinates, any $y \in \R^{s+1}$, and $l \in \N_s$ the
following identity holds
$$
\frac{V\big(\substack{y_1,...,y_s \\
x_1,...,x_s} \big)}{V(x_1,...,x_s)} -
\frac{V\big(\substack{y_1,...,y_{l-1}, y_{s+1}, y_{l+1},...,y_s \\
x_1,...,x_{l-1}, x_{s+1}, x_{l+1},...,x_s})}{V(x_1,...,x_{l-1},
x_{s+1}, x_{l+1},...,x_s)} = (x_l - x_{s+1})
\frac{V\big(\substack{y_1,...,y_l, y_{s+1}, y_{l+1},...,y_s \\
x_1,...,x_l, x_{s+1}, x_{l+1},...,x_s})}{V(x_1,...,x_l, x_{s+1},
x_{l+1},...,x_s)}.
$$
\end{lemma}

\begin{proof}
We consider both sides of the above identity as a multivariate
polynomial in the variables $y_1,...,y_s,y_{s+1}$ and show that
the coefficients in front of $y_k$ on both sides are equal for all
$k \in \N_{s+1}$. Notice first that
\begin{align*}
 V(x_1,...,x_{l-1}, x_{s+1}, x_{l+1},...,x_s) &=
(-1)^{s-l} V(x_1,...,x_{l-1}, x_{l+1},...,x_s, x_{s+1}), \\
V\big(\substack{y_1,...,y_{l-1}, y_{s+1}, y_{l+1},...,y_s \\
x_1,...,x_{l-1}, x_{s+1}, x_{l+1},...,x_s}) &= (-1)^{s-l}
V\big(\substack{y_1,...,y_{l-1}, y_{l+1},...,y_s, y_{s+1} \\
x_1,...,x_{l-1}, x_{l+1},...,x_s, x_{s+1}}).
\end{align*}
We consider four cases according to the partition $\N_{s+1} =
\{1,...,l-1\} \cup \{l\} \cup \{l+1,....,s\}  \cup \{s+1\}$. (In
all product formulae below, it will be assumed that the index $j <
i$. This conditions is omitted for typographical reasons. Also a
hat on top of a multiple in a product denotes that the multiple is
missing.) First, let $k \in \{1,...,l-1\}$. The coefficient in
front of $y_k$ in the expression on the left-hand side is equal to
\begin{align*}
(-1)^{k+1} & \frac{\prod_{i,j \in \N_{s+1}\hspace{-0.1cm}
\backslash \{k,s+1\}}(x_j-x_i) }{\prod_{i,j \in \N_{s+1}
\hspace{-0.1cm} \backslash \{s+1\}}(x_j-x_i)} - (-1)^{k+1}
\frac{\prod_{i,j \in \N_{s+1} \hspace{-0.1cm} \backslash
\{k,l\}}(x_j-x_i)}{\prod_{i,j
\in \N_{s+1} \hspace{-0.1cm} \backslash \{l\}}(x_j-x_i)} \\[0.1cm]
&= \frac{(-1)^{k+1}}{(x_1-x_k) \cdots (x_{k-1}-x_k)(x_k-x_{k+1})
\cdots (x_k-x_s)} \\[0.1cm]
& \hspace{4cm} - \frac{(-1)^{k+1}}{(x_1-x_k) \cdots
(x_{k-1}-x_k)(x_k-x_{k+1}) \cdots \widehat{(x_k-x_{l})}
\cdots (x_k-x_{s+1})} \\[0.1cm]
&= \frac{(-1)^{k+1}}{(x_1-x_k) \cdots (x_{k-1}-x_k)(x_k-x_{k+1})
\cdots \widehat{(x_k-x_{l})} \cdots (x_k-x_s)}
\left(\frac{1}{x_k-x_l} -
\frac{1}{x_k-x_{s+1}} \right) \\[0.1cm]
&=\frac{(-1)^{k+1}(x_l-x_{s+1})}{(x_1-x_k) \cdots
(x_{k-1}-x_k)(x_k-x_{k+1}) \cdots (x_k-x_{s+1})} \\[0.1cm]
&=(-1)^{k+1}(x_l-x_{s+1}) \frac{\prod_{i,j \in
\N_{s+1}\hspace{-0.1cm} \backslash \{k\}}(x_j-x_i) }{\prod_{i,j
\in \N_{s+1}}(x_j-x_i)},
\end{align*}
which is the coefficient in front of $y_k$ on the right-hand side
of the identity.

Suppose now, $k=l$. Then, the coefficient of $y_l$ in the
left-hand side of the identity is
\begin{align*}
(-1)^{l+1} \frac{\prod_{i,j \in \N_{s+1}\hspace{-0.1cm} \backslash
\{l,s+1\}}(x_j-x_i) }{\prod_{i,j \in \N_{s+1} \hspace{-0.1cm}
\backslash \{s+1\}}(x_j-x_i)} - 0 &=
\frac{(-1)^{l+1}}{(x_1-x_l) \cdots (x_{l-1}-x_l)(x_l-x_{l+1}) \cdots (x_l-x_{s})} \\[0.1cm]
&= \frac{(-1)^{l+1}(x_l-x_{s+1})}{(x_1-x_l) \cdots (x_{l-1}-x_l)(x_l-x_{l+1}) \cdots (x_l-x_{s+1})} \\[0.1cm]
&=(-1)^{l+1}(x_l-x_{s+1}) \frac{\prod_{i,j \in
\N_{s+1}\hspace{-0.1cm} \backslash \{l\}}(x_j-x_i) }{\prod_{i,j
\in \N_{s+1}}(x_j-x_i)}.
\end{align*}
When $k \in \{l+1,...,s\}$, the coefficient in front of $y_k$ in
the left-hand side of the identity is:
\begin{align*}
(-1)^{k+1} & \frac{\prod_{i,j \in \N_{s+1}\hspace{-0.1cm}
\backslash \{k,s+1\}}(x_j-x_i) }{\prod_{i,j \in \N_{s+1}
\hspace{-0.1cm} \backslash \{s+1\}}(x_j-x_i)} - (-1)^{k+2}
\frac{\prod_{i,j \in \N_{s+1} \hspace{-0.1cm} \backslash
\{k,l\}}(x_j-x_i)}{\prod_{i,j
\in \N_{s+1} \hspace{-0.1cm} \backslash \{l\}}(x_j-x_i)} \\[0.1cm]
&= \frac{(-1)^{k+1}}{(x_1-x_k) \cdots (x_{k-1}-x_k)(x_k-x_{k+1})
\cdots (x_k-x_s)} \\[0.1cm]
& \hspace{4cm} - \frac{(-1)^{k+2}}{(x_1-x_k) \cdots
\widehat{(x_l-x_{k})} \cdots
(x_{k-1}-x_k)(x_k-x_{k+1}) \cdots (x_k-x_{s+1})} \\[0.1cm]
&= \frac{(-1)^{k+1}}{(x_1-x_k) \cdots \widehat{(x_l-x_k)} \cdots
(x_{k-1}-x_k)(x_k-x_{k+1}) \cdots (x_k-x_s)}
\left(\frac{1}{x_l-x_k} +
\frac{1}{x_k-x_{s+1}} \right) \\[0.1cm]
&=\frac{(-1)^{k+1}(x_l-x_{s+1})}{(x_1-x_k) \cdots
(x_{k-1}-x_k)(x_k-x_{k+1}) \cdots (x_k-x_{s+1})} \\[0.1cm]
&=(-1)^{k+1}(x_l-x_{s+1}) \frac{\prod_{i,j \in
\N_{s+1}\hspace{-0.1cm} \backslash \{k\}}(x_j-x_i) }{\prod_{i,j
\in \N_{s+1}}(x_j-x_i)},
\end{align*}
which is the coefficient in front of $y_k$ on the right-hand side
of the identity.

Finally, when $k=s+1$ the coefficient of $y_{s+1}$ in the
left-hand side of the identity is
\begin{align*}
0-(-1)^{l+1}(-1)^{s-l} \frac{\prod_{i,j \in
\N_{s+1}\hspace{-0.1cm} \backslash \{l,s+1\}}(x_j-x_i)
}{\prod_{i,j \in \N_{s+1} \hspace{-0.1cm} \backslash
\{l\}}(x_j-x_i)} &= \frac{(-1)^{s+2}}{(x_1-x_{s+1}) \cdots
\widehat{(x_l-x_{s+1})} \cdots (x_{s}-x_{s+1})} \\[0.1cm]
&= \frac{(-1)^{s+2}(x_l-x_{s+1})}{(x_1-x_{s+1}) \cdots (x_s-x_{s+1})} \\[0.1cm]
&=(-1)^{s+2}(x_l-x_{s+1}) \frac{\prod_{i,j \in
\N_{s+1}\hspace{-0.1cm} \backslash \{s+1\}}(x_j-x_i) }{\prod_{i,j
\in \N_{s+1}}(x_j-x_i)},
\end{align*}
which is again the coefficient of $y_{s+1}$ on the right. \hfill
\qed
\end{proof}

\begin{theorem}
\label{sym-div-dif} Suppose $g \in C^{k}(I)$. Then, for every
permutation $\sigma \in \tilde{P}^s$, $2 \le s \le k+1$, and every
vector $(x_1,...,x_s)$ with distinct coordinates, we have the
formula
\begin{equation}
\label{Vand-rep} g^{[\sigma]} (x_1,...,x_s) =
\frac{V\big(\substack{g(x_1),...,g(x_s) \\
x_1,...,x_s }\big)}{V(x_1,...,x_s)}.
\end{equation}
In particular, $g^{[\sigma]} (x_1,...,x_s)$ is symmetric
everywhere in its domain.
\end{theorem}

\begin{proof}
The proof is by induction on $s$. When $s=2$ and $x_1 \not= x_2$,
then by the definition we have the representation
$$
g^{[(12)]}(x_1,x_2)= \frac{\left| \begin{array}{cc} g(x_1) &
g(x_2) \\ 1 & 1 \end{array} \right|}{\left| \begin{array}{cc} x_1
& x_2 \\ 1 & 1 \end{array} \right|} = \frac{V(g(x_1),g(x_2);
x_1,x_2)}{V(x_1,x_2)}.
$$
Suppose Representation~(\ref{Vand-rep}) holds for $s$, $2 \le s <
k+1$. Fix a permutation $\sigma \in \tilde{P}^s$ and an $l \in
\N_s$, then $\sigma_{(l)} \in \tilde{P}^{s+1}$. Let
$y=(g(x_1),...,g(x_s),g(x_{s+1}))$. Using
Definition~(\ref{ind-fn-defn}) for any point
$(x_1,...,x_s,x_{s+1})$ with distinct coordinates together with
Lemma~\ref{vand-ident} and the induction hypothesis, we get
\begin{align*}
g^{[\sigma_{(l)}]} &(x_1,...,x_s,x_{s+1}) = \frac{g^{[\sigma]}
(x_1,...,x_s) - g^{[\sigma]}
(x_1,...,x_{l-1},x_{s+1}, x_{l+1},...,x_{s})}{x_l-x_{s+1}} \\[0.1cm]
&= \frac{1}{(x_l-x_{s+1})} \left(
\frac{V\big(\substack{y_1,...,y_s
\\ x_1,...,x_s}\big)}{V(x_1,...,x_s)} -
\frac{V\big(\substack{y_1,...,y_{l-1},y_{s+1}, y_{l+1},...,y_{s} \\
x_1,...,x_{l-1},x_{s+1}, x_{l+1},...,x_{s}}\big)}
{V(x_1,...,x_{l-1},x_{s+1}, x_{l+1},...,x_{s})} \right) \\[0.1cm]
&= \frac{V\big(\substack{y_1,...,y_l, y_{s+1}, y_{l+1},...,y_s \\
x_1,...,x_l, x_{s+1}, x_{l+1},...,x_s})}{V(x_1,...,x_l, x_{s+1},
x_{l+1},...,x_s)} \\
&= \frac{V\big(\substack{y_1,..., y_{s+1} \\
x_1,..., x_{s+1}})}{V(x_1,...,x_{s+1})}.
\end{align*}
Since $\tilde{P}^{s+1} = \{\sigma_{(l)} \, | \, \sigma \in
\tilde{P}^s, l \in \N_s \}$ the induction step is completed.
Finally, using the continuity of $g^{[\sigma]} (x_1,...,x_s)$
shows that it is a symmetric function everywhere on its domain.
\hfill \qed
\end{proof}

Now, we simplify Theorem~\ref{sep-main} significantly. Define the
$(k+1)$-tensor valued map
\begin{align}
\mathbf{g} : \R^n & \,\rightarrow T^{k+1,n}, \mbox{ where } \nonumber \\[-0.3cm]
\label{ILJen-2}  \\[-0.3cm]
\big(\mathbf{g}(\mu)\big)^{i_1...i_{k+1}} &:=
\frac{V\big(\substack{g(\mu_{i_1}),...,g(\mu_{i_{k+1}}) \\
\mu_{i_1},...,\mu_{i_{k+1}}
}\big)}{V(\mu_{i_1},...,\mu_{i_{k+1}})}. \nonumber
\end{align}
Technically, this definition is good only at $\mu$'s with distinct
coordinates, but Lemma~\ref{sym-div-dif} shows that it can be
extended continuously everywhere. Clearly, if $(i_1,...,i_{k+1})
\sim_{\mu} (j_1,...,j_{k+1})$, then
$\big(\mathbf{g}(\mu)\big)^{i_1...i_{k+1}}
=\big(\mathbf{g}(\mu)\big)^{j_1...j_{k+1}}$, which shows that
$\mathbf{g}(\mu)$ is a $\mu$-block-constant tensor for every
$\mu$. Moreover, $\mathbf{g}(\mu)$ is a symmetric tensor, and the
map $\mathbf{g}:\R^n \rightarrow T^{k+1,n}$ is continuous.

\begin{theorem}
\label{sep-main-2} Let $g$ be a $C^k$ function defined on an
interval $I$. Let $X$ be a symmetric matrix with eigenvalues in
the interval $I$, and let $V$ be an orthogonal matrix such that
$X=V(\Diag \lambda(X))V^T$. Then, the matrix valued function $G$
defined by (\ref{G-defn-1}) and (\ref{G-defn-2}) is $C^k$ at $X$.
We have the formula
\begin{equation}
\label{sep-main-form-2} \nabla^{k} G(X) = V\Big(\hspace{0cm}
\sum_{\substack{\sigma \in \tilde{P}^{k+1}}} \hspace{0cm}
\Diag^{\sigma} \mathbf{g}(\lambda(X)) \Big)V^T,
\end{equation}
where the $(k+1)$-tensor valued maps $\mathbf{g}(\cdot)$ is
defined by Equation~(\ref{ILJen-2}).
\end{theorem}

The next corollary is a generalization of Formula~(V.22) from
\cite{Bhatia:1997}.  It is a specialization of the last theorem to
the case when $k=2$. Since $G$ is a symmetric matrix valued
function, the second derivative $\nabla^{2} G(\Diag \mu)[H_1,H_2]$
can be viewed as a symmetric matrix. For every $i=1,...,n$, define
the projection onto the $i$-th coordinate axis
\begin{align*}
P_i : \R^n &\rightarrow \R^n \\
P_i(x) &= x_ie^i.
\end{align*}

\begin{corollary} For $g \in C^2(I)$ and any $n \times n$ symmetric matrices
$H_1$, $H_2$, $H_3$ we have
\begin{align*}
\langle \nabla^{2} G(X)[H_1,H_2], H_3 \rangle &=
2\sum_{p_1,p_2,p_3=1}^{n,n,n}
\mathbf{g}(\lambda(X))^{p_1p_2p_3}\tilde{H}_1^{p_1p_3}
\tilde{H}_2^{p_2p_1}\tilde{H}_3^{p_3p_2}, \\
\nabla^{2} G(X)[H_1,H_2] &= 2\sum_{p_1,p_2,p_3=1}^{n,n,n}
\mathbf{g}(\lambda(X))^{p_1p_2p_3}P_{p_1}\tilde{H}_1P_{p_2}
\tilde{H}_2P_{p_3},
\end{align*}
where $X=V(\Diag \lambda(X))V^T$, and $\tilde{H}_i=V^TH_iV$,
$i=1,2,3$.
\end{corollary}

\begin{proof}
Suppose first that $X=\Diag \mu$ for some $\mu \in
\R^n_{\downarrow}$.
\begin{align*}
\langle \nabla^{2} G(\Diag \mu)[H_1,H_2], H_3 \rangle &=
\nabla^{2} G(\Diag \mu)[H_1,H_2, H_3] \\
&= \Big(\hspace{0cm}
\sum_{\substack{\sigma \in \tilde{P}^{3}}} \hspace{0cm}
\Diag^{\sigma} \mathbf{g}(\mu) \Big)[H_1,H_2, H_3] \\
&= \sum_{\substack{\sigma \in \tilde{P}^{3}}} \langle
\mathbf{g}(\mu), H_1 \circ_{\sigma} H_2 \circ_{\sigma} H_3 \rangle\\
&= \langle \mathbf{g}(\mu), H_1 \circ_{(123)} H_2 \circ_{(123)}
H_3 \rangle + \langle \mathbf{g}(\mu), H_1 \circ_{(132)} H_2
\circ_{(132)} H_3 \rangle \\
&= \hspace{-0.2cm} \sum_{p_1,p_2,p_3=1}^{n,n,n}
\mathbf{g}(\mu)^{p_1p_2p_3}H_1^{p_1p_3}H_2^{p_2p_1}H_3^{p_3p_2} +
\hspace{-0.2cm} \sum_{q_1,q_2,q_3=1}^{n,n,n}
\mathbf{g}(\mu)^{q_1q_2q_3}H_1^{q_1q_2}H_2^{q_2q_3}H_3^{q_3q_1}.
\end{align*}
After re-parametrization of the second sum ($p_1=q_2$, $p_2=q_3$,
$p_3=q_1$), and using the fact that $\mathbf{g}(\mu)$ is a
symmetric tensor, we continue
\begin{align*}
&=\sum_{p_1,p_2,p_3=1}^{n,n,n}
(\mathbf{g}(\mu)^{p_1p_2p_3}+\mathbf{g}(\mu)^{p_3p_1p_2})H_1^{p_1p_3}H_2^{p_2p_1}H_3^{p_3p_2}
= 2\sum_{p_1,p_2,p_3=1}^{n,n,n}
\mathbf{g}(\mu)^{p_1p_2p_3}H_1^{p_1p_3}H_2^{p_2p_1}H_3^{p_3p_2}.
\end{align*}
To show the second representation of $\nabla^{2} G(\Diag
\mu)[H_1,H_2]$ and the general case, when $X$ is not an ordered
diagonal matrix, is routine. \hfill \qed
\end{proof}

\section{The Hessian of spectral functions, revisited}
\label{applic-sect}

In this last section, we illustrate one more time the machinery
developed so far. We recalculate the Hessian of a general spectral
functions at an arbitrary matrix.

One of the strengths of the new approach is that one doesn't need
to have a preconceived notion about the form of the these
derivatives. (Recall that in \cite{LewisSendov:2000a} the formula
for the Hessian of the spectral function was first stated and,
then it was proven that is indeed the correct one. The hind sight
for that formula came from \cite{LewisSendov:2000}.) Here, we
simply differentiate applying the rules developed so far to arrive
at the correct formula. The approach also clearly shows where the
different pieces of the Hessian come from. This should make the
calculation routine and more clear.

\subsection{Two matrix valued maps}
\label{diff-symm:sec}

Let $\mu \in \R^n \rightarrow T(\mu) \in T^{1,n}$, be a
$\mu$-symmetric, differentiable, $1$-tensor-valued map. (In the
next section, $T(\mu) = \nabla f(\mu)$, where $f$ is asymmetric
$C^2$ function.) We define two matrix valued maps $D_0T$ and
$D_1T$ that play an important role in the description of the
Hessian of spectral functions. First
$$
D_0T(\mu) = \nabla T(\mu),
$$
or in other words
$$
D_0T(\mu)^{i_1i_2}= \frac{\partial}{\partial
\mu_{i_2}}(T(\mu)^{i_1}).
$$
Next, define the matrix $D_1T(\mu)$ as follows
$$
D_1T(\mu)^{i_1i_2} = \left\{
\begin{array}{ll}
0, & \mbox{ if } i_1 = i_{2}, \\[0.3cm]
(\nabla T)^{i_1i_1}(\mu) - (\nabla
T)^{i_1i_2}(\mu), & \mbox{ if } i_1 \sim i_{2}, \\[0.3cm]
\displaystyle\frac{T^{i_2}(\mu) - T^{i_1}(\mu)}{\mu_{i_{2}} -
\mu_{i_1}}, & \mbox{ if } i_1 \not\sim i_{2},
\end{array}
\right.
$$
where the equivalence relation is with respect to the vector
$\mu$.  Several of the properties of $D_1T$ are easily seen from
the following integral representations.

\begin{lemma}
If $T(\mu) \in T^{1,n}$ is continuously differentiable, and
$\mu$-symmetric map, then for every $i_1$, $i_2$ $\in \{1,...,n
\}$ we have the representation
\begin{align*}
D_1T(\mu)^{i_1i_2} &= \int_0^1 (\nabla T)^{i_1i_1}(\cdots,
\mu_{i_1}+t(\mu_{i_2}-\mu_{i_1}),\cdots,\mu_{i_2}+t(\mu_{i_1}-\mu_{i_2}),\cdots) - \\
& \hspace{4.7cm} (\nabla T)^{i_1i_2} (\cdots,
\mu_{i_1}+t(\mu_{i_2}-\mu_{i_1}),\cdots,
\mu_{i_2}+t(\mu_{i_1}-\mu_{i_2}),\cdots)dt,
\end{align*}
where the first displayed argument is in position $i_1$ and the
second displayed argument is in position $i_2$. The missing
arguments are the corresponding, unchanged, entries of $\mu$.
\end{lemma}

\begin{proof}
The first case, when $i_1=i_2$ is immediate. In the second, $i_l
\sim i_{2}$ implies that $\mu_{i_{1}} = \mu_{i_2}$ and the
integrand doesn't depend on $t$. In the third case, $i_1 \not\sim
i_{2}$, we can compute the integral using the Fundamental Theorem
of Calculus:
\begin{align*}
D_lT(\mu)^{i_1i_2} &= \frac{1}{\mu_{i_{2}} - \mu_{i_1}} \int_{0}^1
\frac{\partial}{\partial t} T^{i_1}
\big(...,\mu_{i_1}+t(\mu_{i_2}-\mu_{i_1}),...,\mu_{i_{2}}+t(\mu_{i_{1}}-\mu_{i_{2}}),...\big)
\, dt \\[0.2cm]
&= \frac{T^{i_1} (...,\mu_{i_{2}},...,\mu_{i_{1}},...) - T^{i_1}
\big(...,\mu_{i_1},...,\mu_{i_{2}},...\big)
}{\mu_{i_{2}} - \mu_{i_1}} \\
&= \frac{T^{i_2} (...,\mu_{i_1},...,\mu_{i_2},...) - T^{i_1}
\big(...,\mu_{i_1},...,\mu_{i_2},...\big) }{\mu_{2} - \mu_{i_1}},
\end{align*}
where the last equality follows from the fact that $T(\mu)$ is
$\mu$-symmetric. \hfill \qed
\end{proof}

\begin{lemma}
If $T(\mu)$ is differentiable, then both $D_0T(\mu)$ and
$D_1T(\mu)$ are $\mu$-symmetric maps.
\end{lemma}

\begin{proof}
The fact that $D_0T(\mu)$ is $\mu$-symmetric is
Lemma~\ref{sec1-lastlem}. This implies that if $i_1 \sim j_1$,
then $(\nabla T)^{i_1i_1}(\mu)=(\nabla T)^{j_1j_1}(\mu)$. Also, if
$i_1 \sim j_1$ and $i_2 \sim j_2$ with $i_1 \not= i_2$ and $j_1
\not= j_2$, then $(\nabla T)^{i_1i_2}(\mu)=(\nabla
T)^{j_1j_2}(\mu)$. The fact that $T$ is $\mu$-symmetric implies
that if $i_1 \sim j_1$, then $T^{i_1}(\mu) = T^{j_1}(\mu)$. Now it
is easy to see that $D_1T(\mu)$ is $\mu$-symmetric. \hfill \qed
\end{proof}

We conclude this section with a summary of the properties of
$D_0T(\mu)$ and $D_1T(\mu)$
\begin{itemize}
\item For every $l=0,1$, $D_lT(\mu)$ is a matrix valued,
$\mu$-symmetric map. \\[-0.8cm]
\item For every $l=0,1$, $D_lT(\mu)$ is as smooth as $\nabla
T(\mu)$. In other words, if $\nabla T(\mu)$ is continuous, or
several times (continuously) differentiable, then so is
$D_lT(\mu)$. \\[-0.8cm]
\item In addition, if $T = \nabla f(\mu)$ where $f : \R^n
\rightarrow \R$ is a symmetric $C^2$ function, then for every
$l=0,1$, $D_lT(\mu)$ is a symmetric matrix for every $\mu$.
\end{itemize}

\subsection{$f \circ \lambda$ is twice (continuously)
differentiable if, and only if, $f$ is}

Suppose that $f$ is a symmetric function, twice differentiable at
$\mu \in \R^n_{\downarrow}$. Let $E$ be an arbitrary symmetric
matrix. Using Formula~(\ref{dec14:eqn}) together with
Formula~(\ref{kamen}) we compute:
\begin{align*}
\lim_{m \rightarrow \infty} &\frac{\nabla (f \circ \lambda)(\Diag
\mu + M_m) - \nabla (f \circ \lambda)(\Diag \mu)}{\|M_m\|} \\
& \hspace{2.5cm} = \lim_{m \rightarrow \infty}
\frac{U_m\big(\Diag^{\mbox{\tiny (1)}} \nabla
f(\lambda(\Diag \mu +M_m))\big)U_m^T - \Diag^{\mbox{\tiny (1)}} \nabla f(\mu)}{\|M_m\|} \\
& \hspace{2.5cm} = \lim_{m \rightarrow \infty}
\frac{U_m\big(\Diag^{\mbox{\tiny (1)}} \nabla
f(\mu + h_m + o(\|M_m\|))\big)U_m^T - \Diag^{\mbox{\tiny (1)}} \nabla f(\mu)}{\|M_m\|} \\
& \hspace{2.5cm} = \lim_{m \rightarrow \infty}
\frac{U_m\big(\Diag^{\mbox{\tiny (1)}} (\nabla
f(\mu) + \nabla^2f(\mu)[h_m] + o(\|M_m\|))\big)U_m^T - \Diag^{\mbox{\tiny (1)}} \nabla f(\mu)}{\|M_m\|} \\
& \hspace{2.5cm} = \lim_{m \rightarrow \infty}
\frac{U_m\big(\Diag^{\mbox{\tiny (1)}} (\nabla f(\mu))\big)U_m^T -
\Diag^{\mbox{\tiny (1)}} \nabla
f(\mu)}{\|M_m\|}+U\big(\nabla^2f(\mu)[h]\big)U^T.
\end{align*}
For convenience let $T=\nabla f(\mu)$, a block-constant
$1$-tensor.  Using Corollary~\ref{dec14bc} we see that
\begin{align*}
\lim_{m \rightarrow \infty} \frac{U_m\big(\Diag^{\mbox{\tiny (1)}}
T\big)U_m^T - \Diag^{\mbox{\tiny (1)}} T}{\|M_m\|} =
\big(\Diag^{\mbox{\tiny (12)}} T^{\mbox{\tiny (1)}}_{\mbox{\rm
\scriptsize out}}\big)[M].
\end{align*}
Denote $\mathcal{A}_1 = D_0T$, where the operator $D_0$ is defined
in Section~\ref{diff-symm:sec}. Notice that there is a
block-constant vector $b$ such that $\mathcal{A}_1-\Diag b$ is a
block-constant $2$-tensor. Using this notation and
Corollary~\ref{jan11a} we continue:
\begin{align*}
U\big(\nabla^2f(\mu)[h]\big)U^T &= U\big((\mathcal{A}_1-\Diag b +
\Diag b) [h]\big)U^T \\
&= U\big((\mathcal{A}_1-\Diag b) [h]\big)U^T + U\big((\Diag b)
[h]\big)U^T \\
&= \big(\Diag^{\mbox{\tiny (1)(2)}}(\mathcal{A}_1-\Diag b)\big)[M]
+ \big(\Diag^{\mbox{\tiny (12)}} b^{\mbox{\tiny (1)}}_{\mbox{\rm
\scriptsize in}}\big)[M].
\end{align*}
This, shows that $f \circ \lambda$ is twice differentiable.

In order to prove that $f \circ \lambda$ is twice continuously
differentiable we need to reorganize the pieces. Let
$\mathcal{A}_2 = D_1T$, where the operator $D_1$ is defined in
Section~\ref{diff-symm:sec}. Notice that the sum
$\mathcal{A}_1+\mathcal{A}_2$ is block-constant $2$-tensor. This
means that vector $b$ is (can be chosen) such that
$\mathcal{A}_2+\Diag b$ is block-constant, and
$$
\mathcal{A}_2 + \Diag b = T^{\mbox{\tiny (1)}}_{\mbox{\rm
\scriptsize out}} + b^{(1)}_{\mbox{\rm \scriptsize in}}.
$$
Putting everything together we obtain:
\begin{align*}
\nabla^2(f \circ \lambda)(\Diag \mu) &= \Diag^{\mbox{\tiny (12)}}
T^{\mbox{\tiny (1)}}_{\mbox{\rm \scriptsize out}} +
\Diag^{\mbox{\tiny (1)(2)}}(\mathcal{A}_1-\Diag b) +
\Diag^{\mbox{\tiny (12)}} b^{\mbox{\tiny (1)}}_{\mbox{\rm \scriptsize in}} \\
&= \Diag^{\mbox{\tiny (1)(2)}}(\mathcal{A}_1-\Diag b) +
\Diag^{\mbox{\tiny (12)}}(\mathcal{A}_2+\Diag b) \\
&= \Diag^{\mbox{\tiny (1)(2)}}\mathcal{A}_1 +\Diag^{\mbox{\tiny
(12)}}\mathcal{A}_2.
\end{align*}
In the last equality we used the fact that $\Diag^{\mbox{\tiny
(1)(2)}}(\Diag b) = \Diag^{\mbox{\tiny (12)}}(\Diag b),$ which is
very easy to verify. The discussion in
\cite[Section~6]{Sendov2003a} shows that
\begin{equation}
\label{dec14d:eqn} \nabla^2 (f \circ \lambda)(X) =
V\big(\Diag^{\mbox{\tiny (1)(2)}}\mathcal{A}_1(\lambda(X))
+\Diag^{\mbox{\tiny (12)}}\mathcal{A}_2(\lambda(X))\big)V^T,
\end{equation}
where $X=V(\Diag \lambda(X))V^T$.

Moreover, we showed in Section~\ref{diff-symm:sec} that if $f$ is
$C^2$, then both $\mathcal{A}_1$ and $\mathcal{A}_2$ are
continuous. By Proposition~6.2 in \cite{Sendov2003a} it follows
that $\nabla^2 (f \circ \lambda)$ is continuous. That is, $f$ is
$C^2$ if, and only if, $f \circ \lambda$ is.


%

\end{document}